\documentclass[11pt]{amsart}
\usepackage[margin=1.35in]{geometry}
\usepackage{amscd,amsmath,amsxtra,amsthm,amssymb,stmaryrd,xr,mathrsfs,mathtools,enumerate,commath, comment}
\usepackage{stmaryrd}
\usepackage{xcolor}
\usepackage{commath}
\usepackage{comment}
\usepackage{tikz-cd}
\usepackage{float}
\usepackage{longtable} 
\usepackage{pdflscape} 
\usepackage{booktabs}
\usepackage{hyperref}
\usepackage[capitalize,nameinlink,noabbrev]{cleveref}
\usepackage[normalem]{ulem}

\usepackage{multirow}
\usepackage[utf8]{inputenc}

\definecolor{vegasgold}{rgb}{0.77, 0.7, 0.35}
\definecolor{darkgoldenrod}{rgb}{0.72, 0.53, 0.04}
\definecolor{gold(metallic)}{rgb}{0.83, 0.69, 0.22}
\hypersetup{
 colorlinks=true,
 linkcolor=darkgoldenrod,
 filecolor=brown,      
 urlcolor=gold(metallic),
 citecolor=darkgoldenrod,
 pdftitle={An analogue of the Herbrand-Ribet theorem in graph theory},
 }

\usepackage[all,cmtip]{xy}

\usetikzlibrary{shapes.geometric}
\usetikzlibrary{decorations.markings}
\usetikzlibrary{shapes.geometric}
\tikzset{every loop/.style={min distance=10mm,looseness=10}}

\DeclareFontFamily{U}{wncy}{}
\DeclareFontShape{U}{wncy}{m}{n}{<->wncyr10}{}
\DeclareSymbolFont{mcy}{U}{wncy}{m}{n}
\DeclareMathSymbol{\Sh}{\mathord}{mcy}{"58}
\usepackage[T2A,T1]{fontenc}
\usepackage[OT2,T1]{fontenc}

\newtheorem{theorem}{Theorem}[section]
\newtheorem{lemma}[theorem]{Lemma}

\newtheorem{proposition}[theorem]{Proposition}
\newtheorem{corollary}[theorem]{Corollary}
\newtheorem{definition}[theorem]{Definition}

\numberwithin{equation}{section}

\newtheorem{lthm}{Theorem} 

\theoremstyle{remark}
\newtheorem{remark}[theorem]{Remark}

\begin{document}
\title[An analogue of the Herbrand-Ribet theorem in graph theory]{An analogue of the Herbrand-Ribet theorem in graph theory}

\author[D.~Valli\`{e}res]{Daniel Valli\`{e}res}
\address{Daniel Valli\`{e}res\newline Mathematics and Statistics Department, California State University, Chico, CA 95929, USA}
\email{dvallieres@csuchico.edu}

\author[C.~Wilson]{Chase A. Wilson}
\address{Chase A. Wilson\newline Mathematics and Statistics Department, California State University, Chico, CA 95929, USA}
\email{cawilson1@csuchico.edu}

\begin{abstract}
We study an analogue of the Herbrand-Ribet theorem, and its refinement by Mazur and Wiles, in graph theory.  For an odd prime number $p$, we let $\mathbb{F}_{p}$ and $\mathbb{Z}_{p}$ denote the finite field with $p$ elements and the ring of $p$-adic integers, respectively.  We consider Galois covers $Y/X$ of finite graphs with Galois group $\Delta$ isomorphic to $\mathbb{F}_{p}^{\times}$.  Given a $\mathbb{Z}_{p}$-valued character of $\Delta$, we relate the cardinality of the corresponding character component of the $p$-primary subgroup of the degree zero Picard group of $Y$ to the $p$-adic absolute value of the special value at $u=1$ of the corresponding Artin-Ihara $L$-function.
\end{abstract}

\date{\today} 

\maketitle

\tableofcontents

\section{Introduction} \label{Introduction}
Let $p$ be a rational prime.  If $G$ is a finite abelian group, say with additive notation for now, then we let $A$ be its $p$-primary subgroup and $C$ its maximal $p$-elementary abelian quotient.  In other words, $C = G/pG \simeq A/pA$.  One has a short exact sequence 
\begin{equation} \label{primary_vs_elem}
0 \rightarrow pA \rightarrow A \rightarrow C \rightarrow 0 
\end{equation}
of abelian groups.  Moreover if $A \neq 0$, then $pA \subsetneq A$ and therefore one has $p \mid \# G$ if and only if $A \neq 0$ if and only if $C \neq 0$.  It follows that in order to study the $p$-divisibility of a finite abelian group $G$, one can first study whether the maximal $p$-elementary abelian quotient $C$ is trivial or not.  

Assume now that $p$ is odd and let $\mu_{p}$ be the group of $p$th roots of unity in an algebraic closure $\overline{\mathbb{Q}}$ of $\mathbb{Q}$.  Consider the cyclotomic number field $K = \mathbb{Q}(\mu_{p})$.  It is a Galois extension of $\mathbb{Q}$ for which ${\rm Gal}(K/\mathbb{Q}) \simeq \mathbb{F}_{p}^{\times}$, where $\mathbb{F}_{p}$ denotes, as usual, the finite field with $p$ elements.  From now on, we let $\Delta = {\rm Gal}(K/\mathbb{Q})$.  The class number $h(K) = \#{\rm Cl}(K)$ has always been an object of intense study for many reasons, but especially for its intimate connection with Fermat's last theorem.  In particular, mathematicians have been studying the question of whether or not $p$ divides $h(K)$.  An odd prime $p$ is called regular if $p \nmid h(K)$, and it is currently not known if there are infinitely many regular primes.  There is a close connection between this question and Bernoulli numbers.  One way to define the classical Bernoulli numbers $B_{i} \in \mathbb{Q}$ for $i=0,1,2,\ldots$ is by postulating that $B_{0} = 1$, $B_{1} = -1/2$, $B_{i} = 0$ for odd positive integers $i > 1$, and otherwise for even positive integer $i$, via the recursive formula
\begin{equation} \label{Bernoulli}
B_{i} = - \frac{1}{i+1} \sum_{k=0}^{i-1} \binom{i+1}{k}B_{k}. 
\end{equation}
As usual, we let $\mathbb{Z}_{p}$ denote the ring of $p$-adic integers.  It is known, as a consequence of the von Staudt-Clausen theorem, that $B_{i} \in \mathbb{Z}_{p}$ unless $(p-1) \mid i$, in which case $pB_{i} \in \mathbb{Z}_{p}$ (see \cite[Theorem 5.10]{Washington:1997}).  Kummer proved the following theorem (see \cite[Theorem 5.34]{Washington:1997}).
\begin{theorem}[Kummer's criterion]
An odd prime $p$ is regular if and only if it does not divide the numerator of any of the Bernoulli numbers $B_{p-i}$ for $i$ an odd integer satisfying $3 \le i \le p-2$.
\end{theorem}
The first irregular prime is known to be $37$ as can be checked using Kummer's criterion.  Kummer's criterion was refined by Herbrand in \cite{Herbrand:1932}. In order to state what is now known as Herbrand's theorem, we have to introduce some more notation.  Let $C = {\rm Cl}(K)/{\rm Cl}(K)^{p}$ be the maximal $p$-elementary abelian quotient of ${\rm Cl}(K)$ which we view as an $\mathbb{F}_{p}[\Delta]$-module.  Let $\gamma: \Delta \rightarrow \mathbb{F}_{p}^{\times}$ be the natural group isomorphism satisfying $\sigma(\zeta) = \zeta^{\gamma(\sigma)}$ for all $\sigma \in \Delta$ and all $\zeta \in \mu_{p}$.  We view $\gamma$ as an $\mathbb{F}_{p}$-valued character of $\Delta$, and for $i \in \mathbb{Z}/(p-1)\mathbb{Z}$, we let $e_{i}$ be the usual idempotent corresponding to the character $\gamma^{i}$.  We then have a direct sum decomposition
$$C = \bigoplus_{i \in \mathbb{Z}/(p-1)\mathbb{Z}} e_{i} C.$$
Note that by (\ref{primary_vs_elem}), one has $p \mid h(K)$ if and only if $p \mid \#C$ if and only if at least one of the $e_{i}C$ is non-zero.  It is known that $e_{0}C = e_{1}C = 0$ (see \cite[Proposition 6.16]{Washington:1997}).  It is also known that if $e_{j} C \neq 0$ for some even $j$, then $e_{i} C \neq 0$ for an odd $i$ satisfying $i+j \equiv 1 \pmod{p-1}$ (see \cite[Theorem 10.9]{Washington:1997}).  We can now state Herbrand's theorem which is a refinement of Kummer's criterion (see \cite[Theorem 6.17]{Washington:1997} and \cite[Theorem 6.18]{Washington:1997}).
\begin{theorem}[Herbrand]
Let $i$ be an odd integer satisfying $3 \le i \le p-2$.  If $e_{i} C \neq 0$, then $p$ divides the numerator of the Bernoulli number $B_{p-i}$.
\end{theorem}
In \cite{Ribet:1976}, Ribet proved the converse of Herbrand's theorem to obtain what is now known as the Herbrand-Ribet theorem. 
\begin{theorem}[Herbrand-Ribet, first version]
Let $i$ be an odd integer satisfying $3 \le i \le p-2$.  Then, $e_{i} C \neq 0$ if and only if $p$ divides the numerator of $B_{p-i}$.
\end{theorem}

The Herbrand-Ribet theorem can be formulated instead using generalized Bernoulli numbers as we now recall.  Consider the equivariant $L$-function $\theta_{K/\mathbb{Q}}(s)$ and its special value at $s=0$, which by classical results is given by
\begin{equation*} 
\theta_{K/\mathbb{Q}}(0) = \sum_{a=1}^{p-1} \left( \frac{1}{2}- \frac{a}{p}\right) \sigma_{a}^{-1} \in \mathbb{Q}[\Delta], 
\end{equation*}
where $\sigma_{a} \in \Delta$ is defined via $\zeta \mapsto \sigma_{a}(\zeta) = \zeta^{a}$ for all $\zeta \in \mu_{p}$.  Let $\overline{\mathbb{Q}}_{p}$ be an algebraic closure of the field of $p$-adic numbers $\mathbb{Q}_{p}$.  Every $\overline{\mathbb{Q}}_{p}$-valued characters of $\Delta$ actually takes values in $\mathbb{Z}_{p}$, since $\mathbb{Z}_{p}$ contains the $(p-1)$th roots of unity by Hensel's lemma.  Any such $\mathbb{Z}_{p}$-valued character $\psi:\Delta \rightarrow \mathbb{Z}_{p}^{\times}$ induces a unital ring morphism $\psi:\mathbb{Q}[\Delta] \rightarrow \mathbb{Q}_{p}$ which we denote by the same symbol.  For every non-trivial such character $\psi$, one defines the generalized Bernoulli number $B_{1,\psi}$ via the equality
\begin{equation} \label{gen_bernoulli}
B_{1,\psi} = -\psi^{-1}(\theta_{K/\mathbb{Q}}(0)) = \frac{1}{p}\sum_{a=1}^{p-1}a\psi(\sigma_{a}) \in \mathbb{Q}_{p}. 
\end{equation}
Let now $\omega$ be the composition of $\gamma:\Delta \rightarrow \mathbb{F}_{p}^{\times}$ with the usual Teichm\"{u}ller character $\mathbb{F}_{p}^{\times} \hookrightarrow \mathbb{Z}_{p}^{\times}$.  We can view $\omega$ as a $\mathbb{Z}_{p}$-valued character of $\Delta$, and therefore the numbers $B_{1,\omega^{k}} \in \mathbb{Q}_{p}$ make sense for $k=1,2,3,\ldots,p-2$.  It turns out (see \cite[Corollary 5.15]{Washington:1997}) that if $k$ is odd and $k \not \equiv -1 \pmod{p-1}$, then in fact $B_{1,\omega^{k}} \in \mathbb{Z}_{p}$, and thus one can talk about the $p$-divisibility of those numbers.  Using the classical congruence
\begin{equation*} 
B_{1,\omega^{k}} \equiv \frac{B_{k+1}}{k+1} \pmod{p},
\end{equation*}
valid for positive odd integers $k$ satisfying $k \not \equiv -1 \pmod{p-1}$ (see \cite[Corollary 5.15]{Washington:1997}), one has
\begin{equation} \label{congruences}
B_{1,\omega^{-i}} = B_{1,\omega^{p-1-i}} \equiv \frac{B_{p-i}}{p-i} \pmod{p}
\end{equation}
for all odd integer $i$ satisfying $3 \le i \le p-2$.  The congruences (\ref{congruences}) allow one to reformulate the Herbrand-Ribet theorem as follows instead.
\begin{theorem}[Herbrand-Ribet, second version] \label{sec_version}
Let $i$ be an odd integer satisfying $3 \le i \le p-2$.  Then, $e_{i} C \neq 0$ if and only if $p$ divides $B_{1,\omega^{-i}}$.
\end{theorem}
Building on \cite{Mazur:1977} and \cite{Wiles:1980}, and as a consequence of their work \cite{Mazur-Wiles} on the main conjecture in Iwasawa theory for the base field $\mathbb{Q}$, Mazur and Wiles obtained a more refined result that implies the Herbrand-Ribet theorem.  Instead of the maximal $p$-elementary abelian quotient $C$ of ${\rm Cl}(K)$, consider the $p$-primary subgroup $A$ of ${\rm Cl}(K)$, which is now a $\mathbb{Z}_{p}[\Delta]$-module.  For $i=0,1,\ldots,p-2$, let $e_{i} \in \mathbb{Z}_{p}[\Delta]$ be the idempotent in $\mathbb{Z}_{p}[\Delta]$ corresponding to the character $\omega^{i}$.  We have again a direct sum decomposition
$$A = \bigoplus_{i=0}^{p-2} e_{i}A.$$ 
\begin{theorem}[Mazur-Wiles] \label{mw}
Let $i$ be an odd integer satisfying $3 \le i \le p-2$.  Then
$$\# e_{i} A = |B_{1,\omega^{-i}}|_{p}^{-1}, $$
where $|\cdot|_{p}$ denotes the usual $p$-adic absolute value on $\mathbb{Z}_{p}$.
\end{theorem}

The generalized Bernoulli numbers are intimately related with the special value at $s=0$ of Dirichlet $L$-functions.  This is apparent from the original definition (\ref{gen_bernoulli}) above.  Indeed, from the definition of the equivariant $L$-function, if $\psi$ is a non-trivial $\mathbb{C}$-valued character of $\Delta$, one has
$$L(0,\psi) = \psi^{-1}(\theta_{K/\mathbb{Q}}(0)) = -B_{1,\psi} \in \mathbb{C},$$
where $L(s,\psi)$ denotes the corresponding primitive Dirichlet $L$-function.  From this point of view, the Herbrand-Ribet theorem, and its refinement by Mazur and Wiles, can be thought of as a connection between the character components of $C$, or of $A$, and the $p$-divisibility property of a $p$-adic version of the special value at $s=0$ of Dirichlet $L$-functions.

Now, several of those ingredients have an analogue in graph theory.  For instance, Kummer's criterion is a consequence of the analytic class number formula and its factorization in Galois extensions (see \cite[Theorem 5.34]{Washington:1997}), and there is an analogue of the analytic class number formula (\cite{Hashimoto:1990}) and its factorization for Galois covers of graphs (\cite[$\S$3.2.4]{Pengo/Vallieres:2025}).  Moreover, Herbrand's theorem follows from a theorem of Stickelberger (see \cite[Section 6.3]{Washington:1997}), and in \cite{HMSV:2024}, an analogue of Stickelberger's theorem was obtained in the context of graph theory.  (We note in passing that Stickelberger's theorem is a special case of the Brumer-Stark conjecture which was proved recently by Dasgupta and Kakde and their collaborators.  See \cite{Dasgupta/Kakde:2023} and \cite{Dasgupta/Kakde:2023a}.)  The original approach of Ribet and of Mazur-Wiles relies heavily on modular constructions, but once the main conjecture of Iwasawa theory for the base field $\mathbb{Q}$ is known to be true, \cref{mw} is a consequence of a calculation involving Fitting ideals, a tool that is purely algebraic and also available in graph theory.  Although not necessary for our purpose in this paper, we note in passing that an analogue of Iwasawa theory has emerged recently as well in graph theory  (see for instance \cite{Adachi/Mizuno/Tateno:2024}, \cite{Gambheera/Vallieres:2024}, \cite{Ghosh/Ray:2025}, \cite{Gonet:2022}, \cite{Kleine/Muller:2023}, \cite{Kleine/Muller:2025}, \cite{Kundu/Muller:2024}, \cite{lei2022non}, \cite{mcgownvallieresIII}, and \cite{Ray/Vallieres:2023}).  For analogues of Brauer-Kuroda relations in the context of graph theory, see \cite{Mizuno:2025}, and for connections with isogeny graphs, see \cite{Lei/Muller2}, \cite{Lei/Muller1}, and \cite{Lei/Muller3}.

Thus our goal in this paper is to study an analogue of the Herbrand-Ribet theorem and its refinement by Mazur and Wiles for finite graphs.  We will study Galois covers $Y/X$ of finite graphs for which ${\rm Gal}(Y/X) \simeq \mathbb{F}_{p}^{\times}$, where $p$ is an odd rational prime.  We consider the degree zero Picard group ${\rm Pic}^{0}(Y)$ of $Y$, which is known to be a finite abelian group whose cardinality is the number of spanning trees $\kappa(Y)$ of $Y$.  Just as in the number field situation, we consider the $\mathbb{F}_{p}$-vector space $C = {\rm Pic}^{0}(Y)/p{\rm Pic}^{0}(Y)$.  It is acted upon by the Galois group $\Delta = {\rm Gal}(Y/X)$, and thus becomes an $\mathbb{F}_{p}[\Delta]$-module.  As such, it can also be decomposed into character components so that we have a direct sum decomposition
$$C = \bigoplus_{\psi} e_{\psi}C, $$
where the sum is over all $\mathbb{F}_{p}$-valued characters of $\Delta$, and where $e_{\psi}$ are the primitive orthogonal idempotents in $\mathbb{F}_{p}[\Delta]$ corresponding to such characters $\psi$.  On the other hand, we consider the equivariant Ihara zeta function $\theta_{Y/X}(u) \in 1 + u\mathbb{Z}[\Delta]\llbracket u \rrbracket$ of the Galois cover $Y/X$.  It is known that if $\chi(X)$ denotes the Euler characteristic of $X$, then $\theta_{Y/X}(u)^{-1} = (1-u^{2})^{-\chi(X)} \cdot \eta_{Y/X}(u)$ for some explicit polynomial $\eta_{Y/X}(u) \in \mathbb{Z}[\Delta][u]$, and we shall be interested in the special value
$$\eta_{Y/X}(1) \in \mathbb{Z}[\Delta].$$  
Every $\mathbb{F}_{p}$-valued character of $\Delta$ induces a unital ring morphism $\psi:\mathbb{Z}[\Delta] \rightarrow \mathbb{F}_{p}$ which we denote by the same symbol $\psi$.  In particular, if $\psi$ is such a character which is non-trivial, then we let
$$h_{Y/X}(1,\psi) = \psi(\eta_{Y/X}(1)) \in \mathbb{F}_{p}. $$
The first main result we obtain is the following theorem.
\begin{lthm}[\cref{main11}] \label{main1}
Let $p$ be an odd rational prime, and let $Y/X$ be a Galois cover of finite graphs with Galois group $\Delta$ isomorphic to $\mathbb{F}_{p}^{\times}$.  Let $C$ be the maximal $p$-elementary abelian quotient of ${\rm Pic}^{0}(Y)$.  Then, given any non-trivial $\mathbb{F}_{p}$-valued character $\psi$ of $\Delta$, one has 
$$e_{\psi}C \neq 0 \text{ if and only if } h_{Y/X}(1,\psi) = 0 \text{ in } \mathbb{F}_{p}.$$
\end{lthm}
In fact, we obtain \cref{main1} as a consequence of the refined \cref{main2} below.  We let $A$ be the $p$-primary subgroup of ${\rm Pic}^{0}(Y)$; just as before, it is a $\mathbb{Z}_{p}[\Delta]$-module.  Moreover, any $\mathbb{Z}_{p}$-valued character $\psi$ of $\Delta$ induces a unital ring morphism $\psi:\mathbb{Z}[\Delta] \rightarrow \mathbb{Z}_{p}$, which we denote by the same symbol, and when $\psi$ is non-trivial, we let
$$h_{Y/X}(1,\psi) = \psi(\eta_{Y/X}(1)) \in \mathbb{Z}_{p}.$$
For every such character, we also consider the corresponding idempotent $e_{\psi} \in \mathbb{Z}_{p}[\Delta]$, and the character component $e_{\psi}A$.
\begin{lthm}[\cref{main} and \cref{main22}] \label{main2}
Let $p$ be an odd rational prime, and let $Y/X$ be a Galois cover of finite graphs with Galois group $\Delta$ isomorphic to $\mathbb{F}_{p}^{\times}$.  Let $A$ be the $p$-primary subgroup of ${\rm Pic}^{0}(Y)$.  Then, given any non-trivial $\mathbb{Z}_{p}$-valued character of $\Delta$, one has 
$$ \# e_{\psi}A = |h_{Y/X}(1,\psi)|_{p}^{-1},$$
where $|\cdot|_{p}$ denotes the usual $p$-adic absolute value on $\mathbb{Z}_{p}$.
\end{lthm}

The paper is organized as follows.  In \cref{Preliminaries}, we gather together some well-known results about idempotents in group rings and Fitting ideals that we will use in the paper.  In \cref{Graph theory}, we introduce our main objects of study in graph theory.  On the group theory side, we have the Picard group of degree zero of a finite connected graph which is discussed in \cref{Picard Group}.  We study the equivariant Ihara zeta function associated to an abelian cover of finite graphs in \cref{equivariantI}.  This function was introduced in \cite{HMSV:2024}, but our notation will be slightly different.  The connection between the special value at $u=1$ of the equivariant Ihara zeta function and the Fitting ideal of the Picard group is explained in \cref{main_section}.  We also prove our main results in \cref{main_section}, and we end this paper with a few examples in \cref{examples}.

\subsection*{Acknowledgments}
The first author would like to thank Rusiru Gambheera for various stimulating discussions related to the matters of this paper.  The first author also acknowledges support from an AMS-Simons Research Enhancement Grant for PUI Faculty.
\section{Preliminaries} \label{Preliminaries}

\subsection{Idempotents} \label{idem}
Our main reference for this section is \cite{Webb:2016}.  Throughout this paper, by a ring $R$, we shall always mean a ring with a multiplicative identity element denoted by $1_{R}$ or more simply by $1$ if the ring $R$ is understood from the context.  Subrings are required to contain the multiplicative identity element, and ring morphisms are required to be compatible with multiplicative identity elements.  The symbol $R^{\times}$ denotes the group consisting of the units of $R$.  By an $R$-module, we shall mean a left $R$-module throughout.

If $R$ is a commutative ring, and $e_{1},\ldots,e_{n} \in R$ are non-trivial orthogonal (meaning $e_{i}e_{j}=0$ when $i \neq j$) idempotents satisfying $1 = e_{1} + \ldots + e_{n}$, then not only is $R_{i} := Re_{i}$ an ideal of $R$, but also a commutative ring with $1_{R_{i}} = e_{i}$.  Moreover, one has an isomorphism $R \simeq R_{1} \times \ldots \times R_{n}$ of rings, and also a direct sum decomposition $R = R_{1} \oplus \ldots \oplus R_{n}$ as $R$-modules.  The natural projection map $\pi_{i}:R \rightarrow R_{i}$ defined via $r \mapsto \pi_{i}(r) = re_{i}$ is a surjective ring morphism for which ${\rm ker}(\pi_{i}) = {\rm Ann}_{R}(e_{i})$.  If $M$ is now an $R$-module, then we let $M_{i}:= e_{i}M = \{e_{i} \cdot m : m \in M \}$.  By definition, $M_{i}$ is an $R$-submodule of $M$, but it is also an $R_{i}$-module satisfying $M_{i} = \{m \in M : e_{i} \cdot m = m \}$.  Moreover, one has a decomposition $M = M_{1} \oplus \ldots \oplus M_{n}$ as $R$-modules.  The $R_{i}$-module $M_{i}$ can be obtained via a tensor product as well:  Standard properties of the tensor product gives an isomorphism
\begin{equation} \label{tensor}
M_{i} \simeq R_{i} \otimes_{R}M
\end{equation}
as $R_{i}$-modules.

If $f:M \rightarrow N$ is a morphism of $R$-modules, then its restriction to $M_{i}$, denoted by $f_{i}$, satisfies $f_{i}(M_{i}) \subseteq N_{i}$, and thus induces a morphism of $R_{i}$-modules $f_{i}:M_{i} \rightarrow N_{i}$.  Therefore, for each $i=1,\ldots,n$, we have a functor $M \mapsto M_{i}$ from the category of $R$-modules to the category of $R_{i}$-modules.  The following result, albeit simple and left to the reader, is quite useful.
\begin{proposition} \label{exactness}
Consider any exact sequence 
$$0 \rightarrow M \stackrel{\alpha}{\longrightarrow} M' \stackrel{\beta}{\longrightarrow} M'' \rightarrow 0 $$
of $R$-modules.  With the same notation as above, given any $i \in \{1,\ldots,n \}$, one has an exact sequence
$$0 \rightarrow M_{i} \stackrel{\alpha_{i}}{\longrightarrow} M_{i}' \stackrel{\beta_{i}}{\longrightarrow} M_{i}'' \rightarrow 0 $$
of $R_{i}$-modules.  In other words, the functors $M \mapsto M_{i}$ are exact.
\end{proposition}
In this paper, the above will be applied to group rings.  

\subsection{Group rings}
Our main reference again for this section is \cite{Webb:2016}.  If $R$ is a commutative ring and $G$ a finite abelian group, then $R[G]$ denotes the group ring of $G$ over $R$, and we let $\widehat{G}(R) = {\rm Hom}_{\mathbb{Z}}(G,R^{\times})$.  Elements in $\widehat{G}(R)$ are called $R$-valued characters, and the trivial character will be denoted by $\psi_{0}$ throughout.  Moreover, if $\psi \in \widehat{G}(R)$, then its contragredient character will be denoted by $\psi^{*}$.  It satisfies $\psi^{*}(\sigma) = \psi(\sigma^{-1})$ for all $\sigma \in G$.  Note that $\psi^{*} = \psi^{-1}$, when $\widehat{G}(R)$ is viewed as an abelian group with the usual pointwise multiplication.  

Any $\psi \in \widehat{G}(R)$ induces by $R$-linearity a unique ring morphism 
\begin{equation} \label{ring_morphism}
\psi: R[G] \rightarrow R
\end{equation}
which we denote by the same symbol $\psi$.  More generally, if $\varphi:R \rightarrow S$ is a ring morphism between commutative rings and one is given a group morphism $f:G \rightarrow S^{\times}$, then just as above, there is a unique ring morphism 
\begin{equation} \label{ring_morphism2}
f:R[G] \rightarrow S,
\end{equation}
which we denote simply by the same symbol $f$ if $\varphi$ is understood from the context.  For clarity, the ring morphism $f$ above in (\ref{ring_morphism2}) is defined via
$$\sum_{\sigma \in G}m_{\sigma} \cdot \sigma \mapsto f \left(\sum_{\sigma \in G}m_{\sigma} \cdot \sigma \right) = \sum_{\sigma \in G}\varphi(m_{\sigma}) f(\sigma).$$

Let now $R$ be an integral domain with field of fractions $F$, and let $\bar{F}$ denote an algebraic closure of $F$.  We assume that $R$ satisfies the following two conditions:
\begin{enumerate}
\item $\# G \cdot 1_{R} \in R^{\times}$,
\item the group of $\# G$th roots of unity in $\bar{F}$ is contained in $R$.
\end{enumerate}
In this case, we have $\# \widehat{G}(R) = \# G$, and one has the usual orthogonality relations:  If $\psi_{1},\psi_{2} \in \widehat{G}(R)$, then
\begin{equation} \label{ortho1}
\frac{1}{\# G} \sum_{\sigma \in G}\psi_{1}(\sigma) \psi_{2}^{*}(\sigma) = 
\begin{cases}
1, & \text{ if } \psi_{1} = \psi_{2};\\
0, & \text{ if } \psi_{1} \neq \psi_{2},
\end{cases}
\end{equation}
and if $\sigma_{1}, \sigma_{2} \in G$, then
\begin{equation} \label{ortho2}
\frac{1}{\# G} \sum_{\psi \in \widehat{G}(R)}\psi(\sigma_{1}) \psi(\sigma_{2}^{-1}) = 
\begin{cases}
1, & \text{ if } \sigma_{1} = \sigma_{2};\\
0, & \text{ if } \sigma_{1} \neq \sigma_{2}.
\end{cases}
\end{equation}
Given $\psi \in \widehat{G}(R)$, we let
$$e_{\psi} = \frac{1}{\#G} \sum_{\sigma \in G} \psi(\sigma) \sigma^{-1} \in R[G]. $$
The orthogonality relation (\ref{ortho1}) implies that the $e_{\psi}$ are orthogonal idempotents in $R[G]$, and the orthogonality relation (\ref{ortho2}) implies that
$$1 = \sum_{\psi \in \widehat{G}(R)}e_{\psi}. $$
We are thus precisely in the setting of \cref{idem}, and any $R[G]$-module $M$ has a direct sum decomposition
$$M = \bigoplus_{\psi \in \widehat{G}(R)} e_{\psi}M. $$
The $e_{\psi}M$ will be called the character components of $M$.  Note that if $\psi \in \widehat{G}(R)$, a simple calculation shows that
\begin{equation} \label{eigen}
\lambda e_{\psi} = \psi(\lambda) e_{\psi}
\end{equation}
for all $\lambda \in R[G]$, where $\psi$ denotes the ring morphism (\ref{ring_morphism}).  It follows that if $m \in e_{\psi}M$, then $\lambda m = \psi(\lambda) m$ for all $\lambda \in R[G]$.  For that reason, the character components are also sometimes referred to as eigenspaces.  Another useful consequence of the orthogonality relation (\ref{ortho1}) is the equality
\begin{equation} \label{ev_psi}
\psi_{1}(e_{\psi_{2}}) = 
\begin{cases}
1, & \text{ if } \psi_{1} = \psi_{2};\\
0, & \text{ if } \psi_{1} \neq \psi_{2},
\end{cases}
\end{equation}
valid for all $\psi_{1}, \psi_{2} \in \widehat{G}(R)$, where in the equation above $\psi_{1}$ is understood as being the ring morphism (\ref{ring_morphism}).

In this paper, we will apply this to the situation where $G \simeq \mathbb{F}_{p}^{\times}$, and $R = \mathbb{F}_{p}$ or $R = \mathbb{Z}_{p}$ (the $(p-1)$th roots of unity are in $\mathbb{F}_{p}$ and in $\mathbb{Z}_{p}$ because of Fermat's little theorem and Hensel's lemma, respectively).  The commutative ring $\mathbb{F}_{p}[\Delta]$ is semisimple, but $\mathbb{Z}_{p}[\Delta]$ is not, since $\mathbb{Z}_{p}$ itself is not a semisimple ring.  On the other hand, it follows from (\ref{eigen}) and (\ref{ev_psi}) that for all $\psi \in \widehat{\Delta}(\mathbb{Z}_{p})$, the ring morphism $\psi:\mathbb{Z}_{p}[\Delta] \rightarrow \mathbb{Z}_{p}$ coming from (\ref{ring_morphism}) induces an isomorphism
\begin{equation} \label{imp_iso}
\psi:\mathbb{Z}_{p}[\Delta]e_{\psi} \stackrel{\simeq}{\longrightarrow} \mathbb{Z}_{p}
\end{equation} 
of rings.

\subsection{Fitting ideals}
Although the use of initial Fitting ideals is not strictly necessary to obtain our main results in this paper, we find it convenient to use them nonetheless.  Our main reference for this section is \cite[\S 1.4 and \S 3.1]{Northcott:1976}.  Let $R$ be a commutative ring.  We let $M_{m \times n}(R)$ denote the collection of $m \times n$ matrices with entries in $R$.  If $A \in M_{m \times n}(R)$, and $k$ is any non-negative integer, then $D_{k}(A)$ denotes the $k$th determinantal ideal of $A$, that is $D_{k}(A)$ is the ideal of $R$ generated by the $k \times k$ minors of $A$ with the convention that $D_{0}(A) = R$  and that $D_{k}(A) = 0$ if $k > {\rm min}(m,n)$.  The Laplace expansions for determinants show that one has
$$R = D_{0}(A) \supseteq D_{1}(A) \supseteq \ldots \supseteq D_{k}(A) \supseteq \ldots $$

Assume now that we are given a finite free presentation of an $R$-module $M$, that is an exact sequence
\begin{equation} \label{ffp}
F_{1} \stackrel{f}{\longrightarrow} F_{0} \rightarrow M \rightarrow 0
\end{equation}
of  $R$-modules, where both $F_{1}$ and $F_{0}$ are free of finite rank over $R$.  Assume the rank of $F_{1}$ and $F_{0}$ over $R$ are $n$ and $m$, respectively.  Let $A \in M_{m \times n}(R)$ be a presentation matrix for $M$, that is the matrix of $f$ for any choice of $R$-bases for $F_{1}$ and $F_{0}$.  For a non-negative integer $k$, the $k$th Fitting ideal is defined to be $D_{m-k}(A)$, and \cite[Theorem 1 on page 58]{Northcott:1976} implies that these ideals do not depend on the choice of the finite free presentation and of the presentation matrix; thus, they are ideals associated to the $R$-module $M$.  Those ideals are denoted by ${\rm Fit}_{R}^{k}(M)$ and satisfy
$${\rm Fit}_{R}^{0}(M) \subseteq {\rm Fit}_{R}^{1}(M) \subseteq \ldots \subseteq {\rm Fit}_{R}^{k}(M) \subseteq \ldots $$
In this paper, we will be mainly concerned with the $0$th Fitting ideal of $M$, sometimes called the initial Fitting ideal of $M$, which we will also denote simply by ${\rm Fit}_{R}(M)$.

If $M$ and $N$ are two isomorphic  $R$-modules, then ${\rm Fit}_{R}^{k}(M) = {\rm Fit}_{R}^{k}(N)$, but the converse is not true in general.  On the other hand, if $R$ is a principal ideal domain, and $M$ is a finitely generated torsion $R$-module, then the knowledge of the Fitting ideals amounts to the knowledge of the invariant factors of $M$ and thus determines the isomorphism class of $M$ thanks to the structure theorem for such $R$-modules.  Indeed, if 
$$M \simeq R/(d_{1}) \times \ldots \times R/(d_{n}) $$
for some $d_{i} \in R$ satisfying $d_{1} \mid d_{2} \mid \ldots \mid d_{n}$, then by choosing the obvious finite free presentation of $M$, one calculates
\begin{equation*}
\begin{aligned}
{\rm Fit}_{R}^{0}(M) &= (d_{1}\cdot \ldots \cdot d_{n}),\\ 
{\rm Fit}_{R}^{1}(M) &= (d_{1}\cdot \ldots \cdot d_{n-1}), \\
 &\,\,\,   \vdots\\
{\rm Fit}_{R}^{n}(M) &= R.
\end{aligned}
\end{equation*} 
This last calculation is often used when $R = \mathbb{Z}$ or $\mathbb{Z}_{p}$, the ring of $p$-adic integers, in which case one has in particular the following useful observation.
\begin{proposition} \label{fit1}
Let $R$ be $\mathbb{Z}$ or $\mathbb{Z}_{p}$, and let $M$ be any  $R$-module with finite cardinality (that is a finite abelian group, or a finite abelian $p$-group).  Then ${\rm Fit}^{0}_{R}(M)$ is the ideal of $R$ generated by the cardinality of $M$.  
\end{proposition}
Moreover, Fitting ideals behave well with respect to base change (which follows directly after applying the right exact functor $S \otimes_{R} \_$ to (\ref{ffp})), and we record this result in the following proposition for future reference.
\begin{proposition} \label{fit2}
Let $R$ be any commutive ring and let $M$ be any  $R$-module for which there is a finite free presentation.  If $\varphi:R \rightarrow S$ is any ring morphism, then
$${\rm Fit}_{S}^{k}(S \otimes_{R} M) = \langle\varphi({\rm Fit}_{R}^{k}(M))\rangle, $$
where $\langle\varphi({\rm Fit}_{R}^{k}(M))\rangle$ denotes the ideal in $S$ generated by $\varphi({\rm Fit}_{R}^{k}(M)) \subseteq S$.
\end{proposition}

\section{Graph theory} \label{Graph theory}

\subsection{The Picard group} \label{Picard Group}
Our notation will be almost the same as in \cite{HMSV:2024}, so we will be brief.  Throughout, we use Serre's formalism for graphs (see \cite{Serre:1977} and \cite{Sunada:2013}) and all graphs that will appear in this paper will be finite.  If $Y$ is a finite graph, we let $V_{Y}$ denote the collection of vertices of $Y$, and $\mathbf{E}_{Y}$ the collection of directed edges of $Y$.  Moreover, the incidence map $\mathbf{E}_{Y} \rightarrow V_{Y} \times V_{Y}$ is denoted by $\varepsilon \mapsto (o(\varepsilon),t(\varepsilon))$ and the inversion map $\mathbf{E}_{Y} \rightarrow \mathbf{E}_{Y}$ by $\varepsilon \mapsto \bar{\varepsilon}$.  Given $w \in V_{Y}$, we let $\mathbf{E}_{Y,w} = \{\varepsilon \in \mathbf{E}_{Y} : o(\varepsilon) = w \}$ and for each $w \in V_{Y}$, we set ${\rm val}_{Y}(w) = \# \mathbf{E}_{Y,w}$.  We let ${\rm Div}(Y)$ denote the free abelian group on $V_{Y}$, and $s:{\rm Div}(Y) \rightarrow \mathbb{Z}$ is the usual augmentation map given by $w \mapsto s(w) = 1$ for $w \in V_{Y}$.  The kernel of $s$ is denoted by ${\rm Div}^{0}(Y)$ so that we have a short exact sequence
\begin{equation} \label{ses_first}
0 \rightarrow {\rm Div}^{0}(Y) \rightarrow {\rm Div}(Y) \rightarrow \mathbb{Z} \rightarrow 0
\end{equation}
of $\mathbb{Z}$-modules.  The degree operator $\mathcal{D}$ and the adjacency operator $\mathcal{A}$ on ${\rm Div}(Y)$ are defined as usual, namely $\mathcal{D}(w) = {\rm val}_{Y}(w) \cdot w$ and
$$\mathcal{A}(w) = \sum_{\varepsilon \in \mathbf{E}_{Y,w}} t(\varepsilon). $$
The Laplacian operator on ${\rm Div}(Y)$ is defined to be $\mathcal{L} = \mathcal{D} - \mathcal{A}$.  One sets ${\rm Pr}(Y) = \mathcal{L}({\rm Div}(Y)) \subseteq {\rm Div}^{0}(Y)$, and one defines the $\mathbb{Z}$-modules
$${\rm Pic}^{0}(Y) = {\rm Div}^{0}(Y)/{\rm Pr}(Y) \text{ and } {\rm Pic}(Y) = {\rm Div}(Y)/{\rm Pr}(Y).$$
The short exact sequence (\ref{ses_first}) induces another short exact sequence
\begin{equation} \label{ses_second}
0 \rightarrow {\rm Pic}^{0}(Y) \rightarrow {\rm Pic}(Y) \rightarrow \mathbb{Z} \rightarrow 0
\end{equation}
of $\mathbb{Z}$-modules.  If $Y$ is assumed not only to be finite, but also to be connected, then one has yet another short exact sequence
$$0 \longrightarrow \mathbb{Z} \sum_{w \in V_{Y}} w \longrightarrow {\rm Div}(Y) \stackrel{\mathcal{L}}{\longrightarrow} {\rm Pr}(Y) \longrightarrow 0$$
of $\mathbb{Z}$-modules, and it is known that ${\rm Pic}^{0}(Y)$ is a finite abelian group whose cardinality is $\kappa(Y)$, the number of spanning trees of $Y$ (see \cite[Theorem 2.3]{Gambheera/Vallieres:2024}).

Throughout this paper, by an abelian cover $Y/X$ of finite graphs, we mean a Galois cover whose Galois group of deck transformations is abelian (see \cite[$\S$3]{HMSV:2024} for the precise definition of Galois covers of graphs).  We note in passing that throughout this paper, when we say $Y/X$ is a Galois cover of graphs, it includes the condition that both $X$ and $Y$ are connected.  Given such a cover $Y/X$, we let $G = {\rm Gal}(Y/X)$ be its Galois group.  In this situation, the short exact sequences (\ref{ses_first}) and (\ref{ses_second}) become short exact sequences of $\mathbb{Z}[G]$-modules, since everything is acted upon by $G$.  Note that by definition, the exact sequence
\begin{equation} \label{f_f_p_z}
{\rm Div}(Y) \stackrel{\mathcal{L}}{\longrightarrow} {\rm Div}(Y) \rightarrow {\rm Pic}(Y) \rightarrow 0 
\end{equation}
constitutes a finite free presentation of ${\rm Pic}(Y)$ over $\mathbb{Z}$.  Since $G$ acts freely on the vertices lying above a fixed vertex $v$ of $X$, the $\mathbb{Z}[G]$-module ${\rm Div}(Y)$ is actually a free $\mathbb{Z}[G]$-module of rank $\# V_{X}$.  Thus, the finite free presentation (\ref{f_f_p_z}) in fact becomes a finite free presentation of ${\rm Pic}(Y)$ over $\mathbb{Z}[G]$.

\subsection{The equivariant Ihara zeta function} \label{equivariantI}
For this section, our main references are \cite{Terras:2011} and \cite{HMSV:2024}.  (See also \cite{Bass:1992} and \cite{Stark/Terras:1996,Stark/Terras:2000,Stark/Terras:2007}.)  If $R$ is a commutative ring, then $R\llbracket u \rrbracket$ denotes the commutative ring of formal power series in the variable $u$ and with coefficients in $R$.  Similarly, $R[u]$ denotes the commutative ring of polynomials in $u$ with coefficients in $R$.  The ring $R[u]$ is a subring of $R\llbracket u \rrbracket$.  Given a finite graph $X$, its Ihara zeta function, introduced in \cite{Ihara:1966}, is defined to be the formal power series 
$$Z_{X}(u) = \exp\left(\sum_{m=1}^{\infty}N_{m} \frac{u^{m}}{m} \right) \in 1 + u \mathbb{Q}\llbracket u \rrbracket, $$
where $N_{m}$ denotes the number of closed reduced paths in $X$ of length $m$.  Its Euler product is given by
$$Z_{X}(u) = \prod_{\mathfrak{c}}(1 - u^{{\rm len}(\mathfrak{c})})^{-1}, $$
where the product is over all primes in $X$, and ${\rm len}(\mathfrak{c})$ denotes the length of a prime $\mathfrak{c}$ in $X$.  This last formula actually shows in particular that $Z_{X}(u) \in 1 + u \mathbb{Z}\llbracket u \rrbracket \subseteq \mathbb{Z}\llbracket u \rrbracket^{\times}$.  Fix a labeling of the vertices of $X$, say $V_{X} = \{v_{1},\ldots,v_{g} \}$.  Ihara's determinant formula \cite[Theorem 2.5]{Terras:2011} shows that
$$Z_{X}(u)^{-1} = (1-u^{2})^{-\chi(X)} \cdot {\rm det}_{\mathbb{Z}}(I - Au + (D - I)u^{2}) \in \mathbb{Z}[u], $$
where $\chi(X)$ is the Euler characteristic of $X$, and $D$ and $A$ are the degree and adjacency matrices of $X$, namely the matrices of the $\mathbb{Z}$-operators $\mathcal{D}$ and $\mathcal{A}$ on ${\rm Div}(X)$ for the ordered $\mathbb{Z}$-basis $(v_{1},\ldots,v_{g})$ of ${\rm Div}(X)$, respectively.  The matrix $I$ is simply the $g \times g$ identity matrix.

Let now $Y/X$ be an abelian cover of finite graphs, and let $G = {\rm Gal}(Y/X)$.  We define
\begin{equation} \label{equiv_def}
\eta_{Y/X}(u) = {\rm det}_{\mathbb{Z}[G]}(\mathcal{I} - \mathcal{A}u + (\mathcal{D} - \mathcal{I})u^{2}) \in 1 + u\mathbb{Z}[G][u] \subseteq \mathbb{Z}[G]\llbracket u \rrbracket^{\times}, 
\end{equation}
where the operators $\mathcal{D}$ and $\mathcal{A}$ are now viewed as $\mathbb{Z}[G]$-operators on the free $\mathbb{Z}[G]$-module ${\rm Div}(Y)$ of rank $g$.  The operator $\mathcal{I}$ is simply the identity map.  The equivariant Ihara zeta function $\theta_{Y/X}(u)$ is defined via
\begin{equation} \label{ihara_def}
\theta_{Y/X}(u)^{-1} = (1-u^{2})^{-\chi(X)} \cdot \eta_{Y/X}(u) \in 1 + u\mathbb{Z}[G]\llbracket u \rrbracket \subseteq \mathbb{Z}[G]\llbracket u \rrbracket^{\times}, 
\end{equation}
and thus, in particular $\theta_{Y/X}(u) \in \mathbb{Z}[G]\llbracket u \rrbracket$.  

Let $R$ be any commutative ring.  For any $\psi \in \widehat{G}(R)$, the ring morphism $\psi: \mathbb{Z}[G] \rightarrow R$ (for the unique ring morphism $\mathbb{Z} \rightarrow R$) of (\ref{ring_morphism2}) induces, by applying $\psi$ to each coefficient of a power series in $\mathbb{Z}[G]\llbracket u \rrbracket$, a ring morphism 
\begin{equation} \label{ring_mor_ps}
\psi:\mathbb{Z}[G]\llbracket u \rrbracket \rightarrow R\llbracket u \rrbracket
\end{equation}
which we also denote by the same symbol $\psi$.  
\begin{definition}
Let $R$ be a commutative ring.  For any $\psi \in \widehat{G}(R)$, we define
$$L_{Y/X}(u,\psi) = \psi(\theta_{Y/X}(u)) \in R \llbracket u \rrbracket \text{ and } h_{Y/X}(u,\psi) = \psi(\eta_{Y/X}(u)) \in R[u], $$
where on the right-hand side of both equations, $\psi$ is understood to be the ring morphism (\ref{ring_mor_ps}).
\end{definition}
Note that by definition, we have $L_{Y/X}(u,\psi)^{-1} = (1-u^{2})^{-\chi(X)} \cdot h_{Y/X}(u,\psi)$ for all $\psi \in \widehat{G}(R)$.  
\begin{remark}
In particular, if $\psi \in \widehat{G}(\mathbb{C})$, then it follows from \cite[Theorem 18.15]{Terras:2011} after comparing our definitions above to \cite[Definition 18.13]{Terras:2011} that $L_{Y/X}(u,\psi)$ is precisely the Artin-Ihara $L$-function of \cite[Definition 18.6]{Terras:2011} associated to the representation $\rho$ of $G$ having $\psi$ as character.
\end{remark}
We will be particularly interested in the equivariant special value
\begin{equation} \label{sp_val}
\eta_{Y/X}(1) = {\rm det}_{\mathbb{Z}[G]}(\mathcal{L}) \in \mathbb{Z}[G]
\end{equation}
and the special values $h_{Y/X}(1,\psi) = \psi(\eta_{Y/X}(1)) \in R$ for $\psi \in \widehat{G}(R)$ for various commutative rings $R$.

\subsubsection{Explicit description of the equivariant Ihara zeta function}
This section is partly based on \cite{HMSV:2024}.  Let $Y/X$ be an abelian cover of finite graphs with Galois group $G$, and fix a labeling of $V_{X}$, say $V_{X} = \{v_{1},\ldots,v_{g} \}$.  For each $i=1,\ldots,g$, we fix a vertex $w_{i} \in V_{Y}$ lying above $v_{i}$.  Since for each $i=1,\ldots,g$, the group $G$ acts freely on the vertices lying above a fixed $v_{i}$, one has
$${\rm Div}(Y) = \bigoplus_{i=1}^{g}\mathbb{Z}[G] \cdot w_{i} $$
and $\mathbb{Z}[G] \simeq \mathbb{Z}[G]\cdot w_{i}$ as $\mathbb{Z}[G]$-modules.  In this section, we calculate the matrices $\mathbf{D}$ and $\mathbf{A}$ of the $\mathbb{Z}[G]$-operators $\mathcal{D}$ and $\mathcal{A}$ with respect to the ordered $\mathbb{Z}[G]$-basis $(w_{1},\ldots,w_{g})$, respectively. 

Since we have $\mathcal{D}(w_{i}) = {\rm val}_{Y}(w_{i}) \cdot w_{i}$ and ${\rm val}_{Y}(w_{i}) = {\rm val}_{X}(v_{i})$, the matrix $\mathbf{D}$ is a diagonal matrix $(\mathbf{d}_{ij})$, where 
\begin{equation} \label{deg}
\mathbf{d}_{ii} = {\rm val}_{X}(v_{i}).
\end{equation}  
In other words, it is the same matrix as the degree matrix (over $\mathbb{Z}$) of the base graph $X$.

Given any two vertices $w,w' \in V_{Y}$, we define
$$a_{w}(w') = \# \{\varepsilon \in \mathbf{E}_{Y} : o(\varepsilon) = w' \text{ and } t(\varepsilon) = w \} \in \mathbb{Z}. $$
For any $\sigma \in G$ and any $w, w' \in V_{Y}$, one has
\begin{equation} \label{gr}
a_{w}(\sigma w') = a_{\sigma^{-1} w}(w').
\end{equation}
For $i=1,\ldots,g$, we define a function $\ell_{i}:{\rm Div}(Y) \rightarrow \mathbb{Z}[G]$ via
$$w \mapsto \ell_{i}(w) = \sum_{\sigma \in G} a_{w_{i}}(\sigma w) \cdot \sigma^{-1}. $$
We leave it to the reader to check that $\ell_{i}$ is a morphism of $\mathbb{Z}[G]$-modules.
\begin{proposition} \label{in_term_basis}
For all $D \in {\rm Div}(Y)$, one has
$$\mathcal{A}(D) = \sum_{i=1}^{g}\ell_{i}(D)\cdot w_{i}. $$
\end{proposition}
\begin{proof}
For a vertex $w \in V_{Y}$, one has
\begin{equation*}
\begin{aligned}
\sum_{i=1}^{g}\ell_{i}(w) \cdot w_{i} &= \sum_{i=1}^{g} \sum_{\sigma \in G}a_{w_{i}}(\sigma w) \cdot \sigma^{-1}w_{i} \\
&= \sum_{i=1}^{g} \sum_{\sigma \in G}a_{\sigma^{-1}w_{i}}(w) \cdot \sigma^{-1}w_{i} \\
&= \sum_{w_{0} \in V_{Y}} a_{w_{0}}(w) \cdot w_{0} \\
&= \mathcal{A}(w).
\end{aligned}
\end{equation*}
\end{proof}
It follows from \cref{in_term_basis} that the matrix $\mathbf{A} = (\mathbf{a}_{ij})$ of the $\mathbb{Z}[G]$-operator $\mathcal{A}$ on ${\rm Div}(Y)$ with respect to the ordered $\mathbb{Z}[G]$-basis $(w_{1},\ldots,w_{g})$ is given by
\begin{equation} \label{adj}
\mathbf{a}_{ij} = \ell_{i}(w_{j}) = \sum_{\sigma \in G} a_{w_{i}}(\sigma w_{j}) \sigma^{-1}. 
\end{equation}

\begin{corollary} \label{exp_spe}
Let $Y/X$ be an abelian cover of finite graphs with Galois group $G$.  With the same notation as above, one has $\eta_{Y/X}(u) = {\rm det}((p_{ij}(u)))$, where
\begin{equation*}
p_{ij}(u) = 
\begin{cases}
1 - \ell_{i}(w_{i})u + ({\rm val}_{X}(v_{i}) - 1)u^{2}, &\text{ if } i =j;\\
-\ell_{i}(w_{j})u, & \text{ if } i \neq j,
\end{cases}
\end{equation*}
for $i,j=1,\ldots,g$.  Moreover, one has $\eta_{Y/X}(1) = {\rm det}(p_{ij}(1))$.
\end{corollary}
\begin{proof}
This follows directly from (\ref{equiv_def}) combined with (\ref{deg}) and (\ref{adj}).
\end{proof}

Since the group $G$ is assumed to be abelian throughout, one has that the map $G \rightarrow G$ defined via $\sigma \mapsto \sigma^{-1}$ is a group isomorphism.  Therefore, it induces an involution on the ring $\mathbb{Z}[G]$ which we denote by 
\begin{equation} \label{invo}
\iota:\mathbb{Z}[G] \rightarrow \mathbb{Z}[G]. 
\end{equation}
The inversion map $\mathbf{E}_{Y} \rightarrow \mathbf{E}_{Y}$ implies that one has
\begin{equation} \label{inv}
a_{w}(w') = a_{w'}(w)
\end{equation}
for all $w,w' \in V_{Y}$.  From this fact, we obtain the following result.
\begin{lemma} \label{transpose}
Let $Y/X$ be an abelian cover of finite graphs with Galois group $G$.  With the same notation as above, the matrix $\mathbf{A}$ satisfies $\mathbf{A}^{t} = \iota(\mathbf{A})$, where $\iota$ is the ring involution on $M_{g}(\mathbb{Z}[G])$ induced by (\ref{invo}) above.
\end{lemma}
\begin{proof}
Using (\ref{inv}) and (\ref{gr}) above, we calculate
\begin{equation*}
\begin{aligned}
\mathbf{A}^{t} &=  \iota \left(\sum_{\sigma \in G}a_{w_{j}}(\sigma w_{i}) \cdot \sigma \right) \\
&= \iota \left(\sum_{\sigma \in G}a_{\sigma w_{i}}(w_{j}) \cdot \sigma \right) \\
&= \iota \left(\sum_{\sigma \in G}a_{w_{i}}(\sigma^{-1} w_{j}) \cdot \sigma \right) \\
&= \iota(\mathbf{A}).
\end{aligned}
\end{equation*}
\end{proof}
The ring involution $\iota$ defined above in (\ref{invo}) induces yet another ring involution
\begin{equation} \label{pol_inv}
\iota:\mathbb{Z}[G][u] \rightarrow \mathbb{Z}[G][u]
\end{equation}
which we denote by the same symbol.
\begin{proposition} \label{fixed_by_involution}
Let $Y/X$ be an abelian cover of finite graphs with Galois group $G$.  With the same notation as above, one has
\begin{equation} \label{un}
\iota\left( \eta_{Y/X}(u)\right) = \eta_{Y/X}(u),
\end{equation}
where $\iota$ is the map defined in (\ref{pol_inv}) above.  Moreover, if $R$ is any commutative ring and $\psi \in \widehat{G}(R)$, then 
\begin{equation} \label{deux}
h_{Y/X}(u,\psi)= h_{Y/X}(u,\psi^{*}).
\end{equation}
\end{proposition}
\begin{proof}
The equality (\ref{un}) follows from \cref{transpose}, whereas (\ref{deux}) follows after noticing that $\psi \circ \iota:\mathbb{Z}[G] \rightarrow R$ is the same ring morphism as the one induced by $\psi^{*}$. 
\end{proof}
Note in particular, that if $\psi \in \widehat{G}(\mathbb{C})$, then $\psi^{*} = \overline{\psi}$, where the bar here means complex conjugation.  Therefore, for each $\psi \in \widehat{G}(\mathbb{C})$, we have
$$\overline{h_{Y/X}(u,\psi)} = h_{Y/X}(u,\overline{\psi}) = h_{Y/X}(u,\psi), $$
and it follows that
$$h_{Y/X}(u,\psi) \in \mathbb{Z}[\psi]^{+}[u] \text{ and } h_{Y/X}(1,\psi) \in \mathbb{Z}[\psi]^{+}, $$
where $\mathbb{Z}[\psi]^{+}$ denotes the ring of integers of the maximal real subfield of the cyclotomic field $\mathbb{Q}(\psi) = \mathbb{Q}(\{\psi(\sigma): \sigma \in G \})$.
\begin{remark}
The special value $\eta_{Y/X}(1) = {\rm det}_{\mathbb{Z}[G]}(\mathcal{L})$ is the same up to a power of $2$ as the special value, denoted by $\theta_{Y/X}^{*}(1)\cdot e$, of \cite[Theorem 4.5]{HMSV:2024}.  Proposition \ref{fixed_by_involution} shows that the equivariant $L$-function, denoted by $\theta_{Y/X}(u)$ and defined at the beginning of \cite[\S 4]{HMSV:2024}, is in fact the inverse of the equivariant Ihara zeta function defined above in (\ref{ihara_def}).
\end{remark}

\section{The analogue of the Herbrand-Ribet theorem} \label{main_section}
Our starting point is the following observation that was pointed out in the paragraph right after \cite[Definition 2.7]{Kataoka:2024}.  Starting with the finite free presentation (\ref{f_f_p_z}) of ${\rm Pic}(Y)$ over $\mathbb{Z}[G]$, the definition of $0$th Fitting ideal gives ${\rm Fit}_{\mathbb{Z}[G]}({\rm Pic}(Y)) = ({\rm det}_{\mathbb{Z}[G]}(\mathcal{L}))$.  But since ${\rm det}_{\mathbb{Z}[G]}(\mathcal{L}) = \eta_{Y/X}(1)$ by (\ref{sp_val}), we have  
\begin{equation} \label{essential}
{\rm Fit}_{\mathbb{Z}[G]}({\rm Pic}(Y)) = (\eta_{Y/X}(1)), 
\end{equation}
which is an equality of ideals in $\mathbb{Z}[G]$.
\begin{remark}
If $R$ is any commutative ring, and $M$ is an $R$-module having a finite free presentation, then it is known that ${\rm Fit}_{R}(M) \subseteq {\rm Ann}_{R}(M)$ (\cite[Theorem 5 on page 60]{Northcott:1976}).  Therefore, a consequence of (\ref{essential}) is that 
$$\eta_{Y/X}(1) \in {\rm Ann}_{\mathbb{Z}[G]}({\rm Pic}(Y)) \subseteq {\rm Ann}_{\mathbb{Z}[G]}({\rm Pic}^{0}(Y)), $$
and (\ref{essential}) implies \cite[Theorem 4.7]{HMSV:2024}.
\end{remark}
Let $p$ be a rational prime and $G$ a finite abelian group.  If $M$ is a $\mathbb{Z}[G]$-module, then we denote the $\mathbb{Z}_{p}[G]$-module $\mathbb{Z}_{p} \otimes_{\mathbb{Z}}M$ simply by $M_{p}$.

\begin{theorem} \label{main}
Let $p$ be an odd prime number, and let $Y/X$ be a Galois cover of finite graphs for which its Galois group $\Delta$ is isomorphic to $\mathbb{F}_{p}^{\times}$.  Let $A$ be the $p$-primary subgroup of ${\rm Pic}^{0}(Y)$ viewed as a $\mathbb{Z}_{p}[\Delta]$-module.  Then, for all non-trivial characters $\psi \in \widehat{\Delta}(\mathbb{Z}_{p})$, one has
$${\rm Fit}_{\mathbb{Z}_{p}}(e_{\psi} A) = (h_{Y/X}(1,\psi)), $$
which is an equality of ideals in $\mathbb{Z}_{p}$.
\end{theorem}
\begin{proof}
Note that since $\mathbb{Z}_{p}$ is flat over $\mathbb{Z}$, it follows that after tensoring with $\mathbb{Z}_{p}$ over $\mathbb{Z}$ the exact sequence (\ref{f_f_p_z}), one obtains another exact sequence
\begin{equation} \label{p_free_presentation}
{\rm Div}_{p}(Y) \stackrel{\mathcal{L}}{\longrightarrow} {\rm Div}_{p}(Y) \rightarrow {\rm Pic}_{p}(Y) \rightarrow 0 
\end{equation}
now of $\mathbb{Z}_{p}[\Delta]$-modules.  The exact sequence (\ref{p_free_presentation}) is a finite free presentation of ${\rm Pic}_{p}(Y)$ over $\mathbb{Z}_{p}[\Delta]$.  Therefore, we have
$$ {\rm Fit}_{\mathbb{Z}_{p}[\Delta]}({\rm Pic}_{p}(Y)) = (\eta_{Y/X}(1)),$$ 
which is an equality of ideals in $\mathbb{Z}_{p}[\Delta]$.  Consider a non-trivial character $\psi \in \widehat{\Delta}(\mathbb{Z}_{p})$.  The base change property of Fitting ideals (\cref{fit2}) for the natural ring morphism $\pi_{\psi}:\mathbb{Z}_{p}[\Delta] \rightarrow \mathbb{Z}_{p}[\Delta]e_{\psi}$ combined with (\ref{tensor}) give
\begin{equation*} 
{\rm Fit}_{\mathbb{Z}_{p}[\Delta]e_{\psi}}(e_{\psi} {\rm Pic}_{p}(Y)) = (\pi_{\psi}(\eta_{Y/X}(1))),
\end{equation*}  
which is an equality of ideals in $\mathbb{Z}_{p}[\Delta]e_{\psi}$.  Using the isomorphism (\ref{imp_iso}) gives
\begin{equation}\label{eq_fitting}
{\rm Fit}_{\mathbb{Z}_{p}}(e_{\psi} {\rm Pic}_{p}(Y)) = (h_{Y/X}(1,\psi)), 
\end{equation}
which is now an equality of ideals in $\mathbb{Z}_{p}$.  After tensoring the short exact sequence (\ref{ses_second}) with $\mathbb{Z}_{p}$ over $\mathbb{Z}$, one gets a short exact sequence
\begin{equation*}
0 \rightarrow A \rightarrow {\rm Pic}_{p}(Y) \rightarrow \mathbb{Z}_{p} \rightarrow 0
\end{equation*}
of $\mathbb{Z}_{p}[\Delta]$-modules, where the action of $\Delta$ on $\mathbb{Z}_{p}$ is the trivial one.  If $\psi \in \widehat{\Delta}(\mathbb{Z}_{p})$ is non-trivial, then \cref{exactness} and the observation that $e_{\psi} \mathbb{Z}_{p} = 0$, since $\psi \neq \psi_{0}$, imply that we get an isomorphism
$$e_{\psi}A \simeq e_{\psi} {\rm Pic}_{p}(Y) $$
of $\mathbb{Z}_{p}$-modules.  Combining this last isomorphism with (\ref{eq_fitting}) above gives the desired result.
\end{proof}
\cref{main2} from \cref{Introduction} is now a direct consequence of \cref{main}.
\begin{corollary}[\cref{main2}] \label{main22}
Let $p$ be an odd rational prime, and let $Y/X$ be a Galois cover of finite graphs with Galois group $\Delta$ isomorphic to $\mathbb{F}_{p}^{\times}$.  Let $A$ be the $p$-primary subgroup of ${\rm Pic}^{0}(Y)$.  Then, given any non-trivial character $\psi \in \widehat{\Delta}(\mathbb{Z}_{p})$, one has 
$$ \# e_{\psi}A = |h_{Y/X}(1,\psi)|_{p}^{-1},$$
where $|\cdot|_{p}$ denotes the usual $p$-adic absolute value on $\mathbb{Z}_{p}$.
\end{corollary}
\begin{proof}
By \cref{fit1}, we have ${\rm Fit}_{\mathbb{Z}_{p}}(e_{\psi}A) = (\# e_{\psi} A)$, and thus \cref{main} implies $(\# e_{\psi} A) = (h_{Y/X}(1,\psi))$, an equality of ideals in $\mathbb{Z}_{p}$.  It follows that 
$$\# e_{\psi} A = |h_{Y/X}(1,\psi)|_{p}^{-1}$$ 
as desired.
\end{proof}
\cref{main1} from \cref{Introduction} follows by reduction modulo $p$ from \cref{main22}.
\begin{theorem}[\cref{main1}] \label{main11}
Let $p$ be an odd rational prime, and let $Y/X$ be a Galois cover of finite graphs with Galois group $\Delta$ isomorphic to $\mathbb{F}_{p}^{\times}$.  Let $C$ be the maximal $p$-elementary abelian quotient of ${\rm Pic}^{0}(Y)$.  Then, given any non-trivial character $\psi \in \widehat{\Delta}(\mathbb{F}_{p})$, one has 
$$e_{\psi}C \neq 0 \text{ if and only if } h_{Y/X}(1,\psi) = 0 \text{ in } \mathbb{F}_{p}.$$
\end{theorem}
\begin{proof}
The natural ring morphism $\pi:\mathbb{Z}_{p} \rightarrow \mathbb{F}_{p}$ induces a group morphism $\pi_{*}:\widehat{\Delta}(\mathbb{Z}_{p}) \rightarrow \widehat{\Delta}(\mathbb{F}_{p})$ via $\psi \mapsto \pi_{*}(\psi) = \pi \circ \psi$.  Because of the Teichm\"{u}ller character $\mathbb{F}_{p}^{\times} \hookrightarrow \mathbb{Z}_{p}^{\times}$, the group morphism $\pi_{*}$ is in fact a group isomorphism.  The ring morphism $\pi$ also induces a ring morphism $\mathbb{Z}_{p}[\Delta] \rightarrow \mathbb{F}_{p}[\Delta]$ which we denote by the same symbol.  Given $\psi \in \widehat{\Delta}(\mathbb{Z}_{p})$, a simple calculation shows that
\begin{equation} \label{idem_comp}
\pi(e_{\psi}) = e_{\pi_{*}(\psi)},
\end{equation}
an equality of idempotents in $\mathbb{F}_{p}[\Delta]$.  Starting with the short exact sequence 
$$0 \rightarrow pA \rightarrow A \rightarrow C \rightarrow 0$$
of $\mathbb{Z}_{p}[\Delta]$-modules coming from (\ref{primary_vs_elem}), given any $\psi \in \widehat{\Delta}(\mathbb{Z}_{p})$, one has yet another short exact sequence
$$0 \rightarrow e_{\psi}pA \rightarrow e_{\psi}A \rightarrow e_{\psi}C \rightarrow 0 $$
because of \cref{exactness}.  Note that by (\ref{idem_comp}), one has $e_{\psi}C = e_{\pi_{*}(\psi)}C$.  By the argument in the first paragraph of this paper, one has $e_{\pi_{*}(\psi)}C \neq 0$ if and only if $e_{\psi}A \neq 0$.  \cref{main22} implies that for non-trivial $\psi \in \widehat{\Delta}(\mathbb{Z}_{p})$, one has $e_{\psi}A \neq 0$ if and only if $|h_{Y/X}(1,\psi)|_{p} \neq 1$ and this last condition happens if and only if $\pi(h_{Y/X}(1,\psi)) = 0$ in $\mathbb{F}_{p}$.  Since
$$h_{Y/X}(1,\pi_{*}(\psi)) = \pi_{*}(\psi)(\eta_{Y/X}(1)) = \pi(\psi(\eta_{Y/X}(1))) = \pi(h_{Y/X}(1,\psi)), $$
we deduce that $e_{\pi_{*}(\psi)}C \neq 0$ if and only if $h_{Y/X}(1,\pi_{*}(\psi)) = 0$ in $\mathbb{F}_{p}$.  To conclude, it suffices to note, as we did at the beginning of the proof, that $\pi_{*}:\widehat{\Delta}(\mathbb{Z}_{p}) \rightarrow \widehat{\Delta}(\mathbb{F}_{p})$ is in fact a group isomorphism.
\end{proof}
For completeness, we indicate what is happening for the trivial character in the following proposition.

\begin{proposition}
Let $p$ be an odd rational prime, and let $Y/X$ be a Galois cover of finite graphs with Galois group $\Delta$ isomorphic to $\mathbb{F}_{p}^{\times}$.  As before, let $A$ be the $p$-primary subgroup of ${\rm Pic}^{0}(Y)$.  For the trivial character $\psi_{0} \in \widehat{\Delta}(\mathbb{Z}_{p})$, one has
$$e_{\psi_{0}}A \simeq {\rm Pic}_{p}^{0}(X) $$
as $\mathbb{Z}_{p}$-modules, where ${\rm Pic}_{p}^{0}(X)$ denotes the $p$-primary subgroup of ${\rm Pic}^{0}(X)$.  In particular, we have $\# e_{\psi_{0}}A = \kappa_{p}(X)$, where $\kappa_{p}(X) = \# {\rm Pic}_{p}^{0}(X)$.
\end{proposition}
\begin{proof}
For any abelian cover $Y/X$ of finite graphs with Galois group $G$, \cite[Proposition 3.5]{Kataoka:2024} implies that
$${\rm Pic}^{0}(X) \simeq N_{G} {\rm Pic}^{0}(Y) $$
as $\mathbb{Z}$-modules, where $N_{G} = \sum_{\sigma \in G}\sigma \in \mathbb{Z}[G]$.  In our case, we have $N_{\Delta} = (p-1) e_{\psi_{0}}$, and since $p-1 \in \mathbb{Z}_{p}^{\times}$, we have isomorphisms
$${\rm Pic}^{0}_{p}(X) \simeq N_{\Delta} A \simeq e_{\psi_{0}}A $$
as we wanted to show. 
\end{proof}

We end this paper with the following observation:  For an abelian cover $Y/X$ of finite graphs with Galois group $\Delta \simeq \mathbb{F}_{p}^{\times}$, one has
$${\rm dim}_{\mathbb{F}_{p}}(C) \ge {\rm dim}_{\mathbb{F}_{p}}({\rm Pic}_{p}^{0}(X)/p{\rm Pic}_{p}^{0}(X)) + \# \{ \psi \in \widehat{\Delta}(\mathbb{F}_{p}) : h_{Y/X}(1,\psi) = 0 \text{ in } \mathbb{F}_{p} \}, $$
but one does not have an equality in general as Example (\ref{example4}) below shows.

\subsection{Examples} \label{examples}
In order to construct Galois covers $Y/X$ of finite graphs with Galois group $\Delta \simeq \mathbb{F}_{p}^{\times}$, we use voltage assignments (see \cite[$\S$3.1.5]{Pengo/Vallieres:2025}).  A voltage assignment on a graph $X = (V_{X},\mathbf{E}_{X})$ with values in a group $G$ consists of a function $\alpha:\mathbf{E}_{X} \rightarrow G$ satisfying $\alpha(\bar{e}) = \alpha(e)^{-1}$ for all $e \in \mathbf{E}_{X}$.  To every such voltage assignment is associated a graph $X(G,\alpha)$ that comes with a natural covering map $X(G,\alpha) \rightarrow X$.  If $X(G,\alpha)$ is connected, then the cover $X(G,\alpha)/X$ is Galois with Galois group isomorphic to $G$.  Any Galois cover $Y/X$ with Galois group $\Delta \simeq \mathbb{F}_{p}^{\times}$ arises (up to isomorphism of covers) from a voltage assignment as explained in \cite[$\S$3]{Sage/Vallieres:2022} for instance, and thus there is no lost in generality in doing so.  To calculate within the $\mathbb{F}_{p}[\Delta]$-module $C$, we use the isomorphism
$$C \simeq {\rm Div}^{0}(Y)/(p{\rm Div}^{0}(Y) + {\rm Pr}(Y))  $$
of $\mathbb{F}_{p}[\Delta]$-modules, and we do the calculations in ${\rm Div}(Y)/(p{\rm Div}^{0}(Y) + {\rm Pr}(Y))$, since
$$ {\rm Div}^{0}(Y)/(p{\rm Div}^{0}(Y) + {\rm Pr}(Y))\le {\rm Div}(Y)/(p{\rm Div}^{0}(Y) + {\rm Pr}(Y)). $$
We find $\varepsilon_{1},\ldots,\varepsilon_{m} \in {\rm Div}^{0}(Y)$ such that 
$$\varepsilon_{1} + M, \ldots, \varepsilon_{m} + M $$
is a $\mathbb{F}_{p}$-basis for ${\rm Div}^{0}(Y)/M \simeq C$, where from now on, we let $M = p{\rm Div}^{0}(Y) + {\rm Pr}(Y)$.  As some elements $\lambda_{i}$ run over $\{0,1,\ldots,p-2 \}$ for $i=1,\ldots,m$, the elements $\alpha + M$, where
$$\alpha = \lambda_{1}\varepsilon_{1} + \ldots + \lambda_{m} \varepsilon_{m}, $$
run over all the elements of $C \simeq {\rm Div}^{0}(Y)/M$.  For each $\psi \in \widehat{\Delta}(\mathbb{F}_{p})$, we calculate a lift $f_{\psi} \in \mathbb{Z}[\Delta]$ of $e_{\psi} \in \mathbb{F}_{p}[\Delta]$.  One has
\begin{equation} \label{isot}
e_{\psi}(\alpha + M) = \alpha + M \text{ if and only if } f_{\psi} \alpha - \alpha \in M. 
\end{equation}
For each non-trivial $\psi \in \widehat{\Delta}(\mathbb{F}_{p})$, we run over all $\alpha$ and we check if condition (\ref{isot}) is satisfied.  All of those calculations can be done with the Smith normal form of a presentation matrix for ${\rm Div}(Y)/M$ and this allows us to calculate $\# e_{\psi} C$.  The special values $\eta_{Y/X}(1)$ and $h_{Y/X}(1,\psi)$ are calculated using \cref{exp_spe}.  All calculations are made using \cite{SAGE}.

\begin{enumerate}
\item Take $X$ to be the bouquet graph on one vertex with two undirected loops, and let $p = 5$.  Let the directed edges of $X$ be $\mathbf{E}_{X} = \{s_{1},s_{2},\bar{s}_{1},\bar{s}_{2} \}$, where $\{s_{1},s_{2} \}$ is a fixed orientation of $X$.  Consider the voltage assignment $\alpha:\mathbf{E}_{X} \rightarrow \mathbb{F}_{5}^{\times}$ defined via $\alpha(s_{1}) = 2$ and $\alpha(s_{2}) = 3$.  The derived graph $Y = X(\mathbb{F}_{5}^{\times},\alpha)$ and the Galois cover $Y \rightarrow X$ can be visualized as follows:  
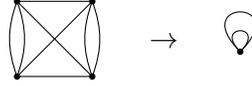
\begin{equation*}  
\hspace{2cm}
\begin{tikzpicture}[baseline={([yshift=-0.7ex] current bounding box.center)}]
\draw[fill=black] (0,0) circle (1pt);
\draw[fill=black] (1,0) circle (1pt);
\draw[fill=black] (0,1) circle (1pt);
\draw[fill=black] (1,1) circle (1pt);
\path (0,0) edge (1,0);
\path (0,1) edge (1,1);
\path (0,1) edge (1,0);
\path (0,0) edge (1,1);
\path (0,0) edge [bend right=20] (0,1);
\path (0,0) edge [bend left=20] (0,1);
\path (1,0) edge [bend right=20] (1,1);
\path (1,0) edge [bend left=20] (1,1);
\end{tikzpicture}
\, \, \, \, \, \, \, \, \rightarrow
\begin{tikzpicture}[baseline={([yshift=-1.7ex] current bounding box.center)}]
\draw[fill=black] (0,0) circle (1pt);
\path (0,0) edge [loop,min distance=15] (0,0);
\path (0,0) edge [loop] (0,0);
\end{tikzpicture}
\end{equation*}
One has ${\rm Pic}^{0}(Y) \simeq \mathbb{Z}/3\mathbb{Z}  \times \mathbb{Z}/12\mathbb{Z}$ as $\mathbb{Z}$-modules so that in particular $C = A = 0$.  Let $\gamma = {\rm id}_{\mathbb{F}_{5}^{\times}}$ viewed as an element of $\widehat{\Delta}(\mathbb{F}_{5})$.  We calculate 
\begin{center}
\begin{tabular}{c|c|c}    
$\psi = \gamma^{i}$ & $e_{\psi}C$ & $h_{Y/X}(1,\psi) \text{ in } \mathbb{F}_{5}$\\
\hline \hline
$1$ & $0$ & $1$ \\
$2$ & $0$ & $4$ \\
$3$ & $0$ & $1$ \\
\end{tabular}
\end{center}
None of the special values is $0$ in $\mathbb{F}_{5}$ as expected by \cref{main11}.  Note also that $\gamma^{*} = \gamma^{3}$ and thus the equality $h_{Y/X}(1,\gamma) = h_{Y/X}(1,\gamma^{3})$ is expected as well by \cref{fixed_by_involution}.

\item Take $X$ to be the graph consisting of two vertices $v_{1}$ and $v_{2}$ with a single undirected loop at $v_{1}$ and three undirected edges between $v_{1}$ and $v_{2}$, and let $p = 5$.  Choose an orientation $\{s_{1},s_{2},s_{3},s_{4} \}$ of $\mathbf{E}_{X}$, where $s_{1}$ is a directed edge corresponding to the single loop, and $s_{i}$ are all going from $v_{1}$ to $v_{2}$ for $i=2,3,4$.  We consider the voltage assignment $\alpha:\mathbf{E}_{X} \rightarrow \mathbb{F}_{5}^{\times}$ defined by
$$\alpha(s_{1}) = \alpha(s_{3}) = 2, \alpha(s_{2}) = 4, \text{ and } \alpha(s_{4}) = 1. $$
The derived graph $Y = X(\mathbb{F}_{5}^{\times},\alpha)$ and the Galois cover $Y \rightarrow X$ can be visualized as follows:
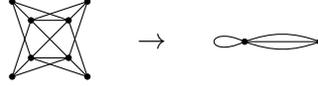
\begin{equation*}  
\begin{tikzpicture}[baseline={([yshift=-0.7ex] current bounding box.center)}]
\draw[fill=black] (0,0) circle (1pt);
\draw[fill=black] (1,0) circle (1pt);
\draw[fill=black] (0,1) circle (1pt);
\draw[fill=black] (1,1) circle (1pt);
\draw[fill=black] (0.25,0.25) circle (1pt);
\draw[fill=black] (0.75,0.25) circle (1pt);
\draw[fill=black] (0.25,0.75) circle (1pt);
\draw[fill=black] (0.75,0.75) circle (1pt);
\path (0.25,0.25) edge (0.75,0.25);
\path (0.25,0.75) edge (0.75,0.75);
\path (0.25,0.25) edge (0.75,0.75);
\path (0.25,0.75) edge (0.75,0.25);
\path (0,0) edge (0.25,0.25);
\path (0,1) edge (0.25,0.75);
\path (1,0) edge (0.75,0.25);
\path (1,1) edge (0.75,0.75);
\path (0,0) edge (0.25,0.75);
\path (0,1) edge (0.75,0.75);
\path (1,0) edge (0.25,0.25);
\path (1,1) edge (0.75,0.25);
\path (0,0) edge (0.75,0.25);
\path (0,1) edge (0.25,0.25);
\path (1,0) edge (0.75,0.75);
\path (1,1) edge (0.25,0.75);
\end{tikzpicture}
\, \, \, \, \, \, \, \, \rightarrow \, \, \, \, \, \, 
\begin{tikzpicture}[baseline={([yshift=-0.5ex] current bounding box.center)}]
\draw[fill=black] (0,0) circle (1pt);
\draw[fill=black] (1,0) circle (1pt);
\path (0,0) edge [loop left, in = 155, out = 205,min distance=6mm] (0,0);
\path (0,0) edge [bend right=20] (1,0);
\path (0,0) edge [bend left=20] (1,0);
\path (0,0) edge (1,0);
\end{tikzpicture}
\end{equation*}
One has $A \simeq  \mathbb{Z}/5\mathbb{Z}$ as $\mathbb{Z}$-modules so that in particular $C = A \simeq \mathbb{Z}/5\mathbb{Z}$.  Let $\gamma = {\rm id}_{\mathbb{F}_{5}^{\times}}$ viewed as an element of $\widehat{\Delta}(\mathbb{F}_{5})$.  We calculate 
\begin{center}
\begin{tabular}{c|c|c}    
$\psi = \gamma^{i}$ & $e_{\psi}C$ & $h_{Y/X}(1,\psi) \text{ in } \mathbb{F}_{5}$\\
\hline \hline
$1$ & $0$ & $4$ \\
$2$ & $\mathbb{F}_{5}$ & $0$ \\
$3$ & $0$ & $4$ \\
\end{tabular}
\end{center}
We have $e_{\psi} \cdot C \neq 0$ if and only if $h_{Y/X}(1,\psi) = 0$ as expected by \cref{main11}.  If $\omega$ is the Teichm\"{u}ller character, then $\omega \circ \gamma^{2}$ is a lift of $\gamma^{2}$, and we calculated
$$h_{Y/X}(1,\omega \circ \gamma^{2}) = 4 \cdot 5 + \ldots \in \mathbb{Z}_{5}, $$
so that $\# e_{\omega \circ \gamma^{2}} A = |h_{Y/X}(1,\omega \circ \gamma^{2})|_{5}^{-1}$ as expected by \cref{main22}. 
\item Take $X$ to be the graph consisting of two vertices $v_{1}$ and $v_{2}$ with a single undirected loop at $v_{1}$, two undirected edges between $v_{1}$ and $v_{2}$, and a single undirected loop at $v_{2}$, and let $p = 11$.  Choose an orientation $\{s_{1},s_{2},s_{3},s_{4} \}$ of $\mathbf{E}_{X}$, where $s_{1}$ is a directed loop at $v_{1}$, and $s_{i}$ are all going from $v_{1}$ to $v_{2}$ for $i=2,3$, and $s_{4}$ is a directed loop at $v_{2}$.  We consider the voltage assignment $\alpha:\mathbf{E}_{X} \rightarrow \mathbb{F}_{11}^{\times}$ defined by
$$\alpha(s_{1}) = 4, \alpha(s_{2}) =  \alpha(s_{3}) = 1, \text{ and } \alpha(s_{4}) = 10. $$
The derived graph $Y = X(\mathbb{F}_{11}^{\times},\alpha)$ and the Galois cover $Y \rightarrow X$ can be visualized as follows:
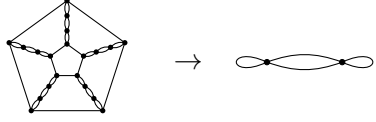
\begin{equation*}
\begin{tikzpicture}[baseline={([yshift=-1.1ex] current bounding box.center)}]
\node[draw=none,minimum size=1.6cm,regular polygon,regular polygon sides=5] (a) {};
\node[draw=none, minimum size=1.2cm,regular polygon,regular polygon sides=5] (b) {};
\node[draw=none, minimum size=0.8cm,regular polygon,regular polygon sides=5] (c) {};
\node[draw=none, minimum size=0.4cm,regular polygon,regular polygon sides=5] (d) {};

\foreach \x in {1,2,...,5}
  \fill (a.corner \x) circle[radius=1pt];
  
\foreach \y in {1,2,...,5}
  \fill (b.corner \y) circle[radius=1pt];

\foreach \y in {1,2,...,5}
  \fill (c.corner \y) circle[radius=1pt];

\foreach \y in {1,2,...,5}
  \fill (d.corner \y) circle[radius=1pt];
  
\path (a.corner 1) edge (a.corner 2);
\path (a.corner 2) edge (a.corner 3);
\path (a.corner 3) edge (a.corner 4);
\path (a.corner 4) edge (a.corner 5);
\path (a.corner 5) edge (a.corner 1);

\path (d.corner 1) edge (d.corner 2);
\path (d.corner 2) edge (d.corner 3);
\path (d.corner 3) edge (d.corner 4);
\path (d.corner 4) edge (d.corner 5);
\path (d.corner 5) edge (d.corner 1);

\path (a.corner 1) edge [bend right=30] (b.corner 1);
\path (a.corner 1) edge [bend left=30] (b.corner 1);
\path (a.corner 2) edge [bend right=30] (b.corner 2);
\path (a.corner 2) edge [bend left=30] (b.corner 2);
\path (a.corner 3) edge [bend right=30] (b.corner 3);
\path (a.corner 3) edge [bend left=30] (b.corner 3);
\path (a.corner 4) edge [bend right=30] (b.corner 4);
\path (a.corner 4) edge [bend left=30] (b.corner 4);
\path (a.corner 5) edge [bend right=30] (b.corner 5);
\path (a.corner 5) edge [bend left=30] (b.corner 5);

\path (c.corner 1) edge [bend right=30] (b.corner 1);
\path (c.corner 1) edge [bend left=30] (b.corner 1);
\path (c.corner 2) edge [bend right=30] (b.corner 2);
\path (c.corner 2) edge [bend left=30] (b.corner 2);
\path (c.corner 3) edge [bend right=30] (b.corner 3);
\path (c.corner 3) edge [bend left=30] (b.corner 3);
\path (c.corner 4) edge [bend right=30] (b.corner 4);
\path (c.corner 4) edge [bend left=30] (b.corner 4);
\path (c.corner 5) edge [bend right=30] (b.corner 5);
\path (c.corner 5) edge [bend left=30] (b.corner 5);

\path (c.corner 1) edge [bend right=30] (d.corner 1);
\path (c.corner 1) edge [bend left=30] (d.corner 1);
\path (c.corner 2) edge [bend right=30] (d.corner 2);
\path (c.corner 2) edge [bend left=30] (d.corner 2);
\path (c.corner 3) edge [bend right=30] (d.corner 3);
\path (c.corner 3) edge [bend left=30] (d.corner 3);
\path (c.corner 4) edge [bend right=30] (d.corner 4);
\path (c.corner 4) edge [bend left=30] (d.corner 4);
\path (c.corner 5) edge [bend right=30] (d.corner 5);
\path (c.corner 5) edge [bend left=30] (d.corner 5);
\end{tikzpicture}
\, \, \, \, \, \, \, \, \rightarrow \, \, \,
\begin{tikzpicture}[baseline={([yshift=-0.6ex] current bounding box.center)}]
\draw[fill=black] (0,0) circle (1pt);
\draw[fill=black] (1,0) circle (1pt);
\path (0,0) edge [loop left, in = 155, out = 205,min distance=6mm] (0,0);
\path (0,0) edge [bend right=20] (1,0);
\path (0,0) edge [bend left=20] (1,0);
\path (1,0) edge [loop right, in = 25, out = 335,min distance=6mm] (1,0);
\end{tikzpicture}
\end{equation*}
One has $A \simeq  \mathbb{Z}/11^{2}\mathbb{Z} \times \mathbb{Z}/11^{2}\mathbb{Z}$ as $\mathbb{Z}$-modules so that in particular $C \simeq \mathbb{Z}/11\mathbb{Z} \times \mathbb{Z}/11\mathbb{Z}$.  Let $\gamma = {\rm id}_{\mathbb{F}_{11}^{\times}}$ viewed as an element of $\widehat{\Delta}(\mathbb{F}_{11})$.  We calculate 
\begin{center}
\begin{tabular}{c|c|c}    
$\psi = \gamma^{i}$ & $e_{\psi}C$ & $h_{Y/X}(1,\psi) \text{ in } \mathbb{F}_{11}$\\
\hline \hline
$1$ & $\mathbb{F}_{11}$ & $0$ \\
$2$ & $0$ & $9$ \\
$3$ & $0$ & $2$ \\
$4$ & $0$ & $1$ \\
$5$ & $0$ & $8$ \\
$6$ & $0$ & $1$ \\
$7$ & $0$ & $2$ \\
$8$ & $0$ & $9$ \\
$9$ & $\mathbb{F}_{11}$ & $0$ \\
\end{tabular}
\end{center}
We have $e_{\psi} \cdot C \neq 0$ if and only if $h_{Y/X}(1,\psi) = 0$ as expected by \cref{main11}.  If $\omega$ is the Teichm\"{u}ller character, then $\omega \circ \gamma$ is a lift of $\gamma$, and we calculated
$$h_{Y/X}(1,\omega \circ \gamma) = 2 \cdot 11^{2} + 7\cdot 11^{3} + 9\cdot 11^{4} + \ldots \in \mathbb{Z}_{11}, $$
so that $11^{2} = \# e_{\omega \circ \gamma} A = |h_{Y/X}(1,\omega \circ \gamma)|_{11}^{-1}$ as expected by \cref{main22}.
\item Take $X$ to be the same graph as in the previous example, and let $p = 11$.  Pick also the same orientation, and consider $\alpha:\mathbf{E}_{X} \rightarrow \mathbb{F}_{11}^{\times}$ defined by
$$\alpha(s_{1}) = \alpha(s_{4})= 2, \alpha(s_{2}) = 1, \text{ and } \alpha(s_{3}) = 10. $$
The derived graph $Y = X(\mathbb{F}_{11}^{\times},\alpha)$ and the Galois cover $Y \rightarrow X$ can be visualized as follows:
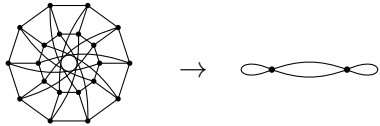
\begin{equation*}
\begin{tikzpicture}[baseline={([yshift=-1.1ex] current bounding box.center)}]
\node[draw=none,minimum size=1.6cm,regular polygon,regular polygon sides=10] (a) {};
\node[draw=none, minimum size=0.8cm,regular polygon,regular polygon sides=10] (b) {};

\foreach \x in {1,2,...,10}
    \fill (a.corner \x) circle[radius=1pt];
    
\foreach \y in {1,2,...,10}
    \fill (b.corner \y) circle[radius=1pt];

\path (a.corner 1) edge (a.corner 2);
\path (a.corner 2) edge (a.corner 3);
\path (a.corner 3) edge (a.corner 4);
\path (a.corner 4) edge (a.corner 5);
\path (a.corner 5) edge (a.corner 6);
\path (a.corner 6) edge (a.corner 7);
\path (a.corner 7) edge (a.corner 8);
\path (a.corner 8) edge (a.corner 9);
\path (a.corner 9) edge (a.corner 10);
\path (a.corner 10) edge (a.corner 1);

\path (b.corner 1) edge (b.corner 2);
\path (b.corner 2) edge (b.corner 3);
\path (b.corner 3) edge (b.corner 4);
\path (b.corner 4) edge (b.corner 5);
\path (b.corner 5) edge (b.corner 6);
\path (b.corner 6) edge (b.corner 7);
\path (b.corner 7) edge (b.corner 8);
\path (b.corner 8) edge (b.corner 9);
\path (b.corner 9) edge (b.corner 10);
\path (b.corner 10) edge (b.corner 1);

\path (a.corner 1) edge (b.corner 1);
\path (a.corner 2) edge (b.corner 2);
\path (a.corner 3) edge (b.corner 3);
\path (a.corner 4) edge (b.corner 4);
\path (a.corner 5) edge (b.corner 5);
\path (a.corner 6) edge (b.corner 6);
\path (a.corner 7) edge (b.corner 7);
\path (a.corner 8) edge (b.corner 8);
\path (a.corner 9) edge (b.corner 9);
\path (a.corner 10) edge (b.corner 10);

\path (a.corner 1) edge [bend right=20] (b.corner 6);
\path (a.corner 2) edge [bend right=20] (b.corner 7);
\path (a.corner 3) edge [bend right=20] (b.corner 8);
\path (a.corner 4) edge [bend right=20] (b.corner 9);
\path (a.corner 5) edge [bend right=20] (b.corner 10);
\path (a.corner 6) edge [bend right=20] (b.corner 1);
\path (a.corner 7) edge [bend right=20] (b.corner 2);
\path (a.corner 8) edge [bend right=20] (b.corner 3);
\path (a.corner 9) edge [bend right=20] (b.corner 4);
\path (a.corner 10) edge [bend right=20] (b.corner 5);

\end{tikzpicture}
\, \, \, \, \, \, \, \, \rightarrow \, \, \,
\begin{tikzpicture}[baseline={([yshift=-0.6ex] current bounding box.center)}]
\draw[fill=black] (0,0) circle (1pt);
\draw[fill=black] (1,0) circle (1pt);
\path (0,0) edge [loop left, in = 155, out = 205,min distance=6mm] (0,0);
\path (0,0) edge [bend right=20] (1,0);
\path (0,0) edge [bend left=20] (1,0);
\path (1,0) edge [loop right, in = 25, out = 335,min distance=6mm] (1,0);
\end{tikzpicture}
\end{equation*}
One has $A \simeq  (\mathbb{Z}/11 \mathbb{Z})^{4}$ as $\mathbb{Z}$-modules so that in particular $A = C \simeq (\mathbb{Z}/11\mathbb{Z})^{4}$.  Let $\gamma = {\rm id}_{\mathbb{F}_{11}^{\times}}$ viewed as an element of $\widehat{\Delta}(\mathbb{F}_{11})$.  We calculate 
\begin{center}
\begin{tabular}{c|c|c}    
$\psi = \gamma^{i}$ & $e_{\psi}C$ & $h_{Y/X}(1,\psi) \text{ in } \mathbb{F}_{11}$\\
\hline \hline
$1$ & $0$ & $5$ \\
$2$ & $0$ & $5$ \\
$3$ & $\mathbb{F}_{11}^{2}$ & $0$ \\
$4$ & $0$ & $8$ \\
$5$ & $0$ & $3$ \\
$6$ & $0$ & $8$ \\
$7$ & $\mathbb{F}_{11}^{2}$ & $0$ \\
$8$ & $0$ & $5$ \\
$9$ & $0$ & $5$ \\
\end{tabular}
\end{center}
We have $e_{\psi} \cdot C \neq 0$ if and only if $h_{Y/X}(1,\psi) = 0$ as expected by \cref{main11}.  If $\omega$ is the Teichm\"{u}ller character, then $\omega \circ \gamma^{3}$ is a lift of $\gamma^{3}$, and we calculated
$$h_{Y/X}(1,\omega \circ \gamma^{3}) = 9 \cdot 11^{2} + 7 \cdot 11^{4} + 9 \cdot 11^{5} \ldots \in \mathbb{Z}_{11}, $$
so that $11^{2} = \# e_{\omega \circ \gamma^{3}} A = |h_{Y/X}(1,\omega \circ \gamma^{3})|_{11}^{-1}$ as expected by \cref{main22}. \label{example4}

\end{enumerate}

\bibliographystyle{plain}
\bibliography{references}

\begin{thebibliography}{10}

\bibitem{Adachi/Mizuno/Tateno:2024}
Taiga Adachi, Kosuke Mizuno, and Sohei Tateno.
\newblock Iwasawa theory for weighted graphs.
\newblock {\em {P}reprint, arXiv:2412.01612}, 2023.

\bibitem{Bass:1992}
Hyman Bass.
\newblock The {I}hara-{S}elberg zeta function of a tree lattice.
\newblock {\em Internat. J. Math.}, 3(6):717--797, 1992.

\bibitem{Dasgupta/Kakde:2023}
Samit Dasgupta and Mahesh Kakde.
\newblock On the {B}rumer-{S}tark conjecture.
\newblock {\em Ann. of Math. (2)}, 197(1):289--388, 2023.

\bibitem{Dasgupta/Kakde:2023a}
Samit Dasgupta, Mahesh Kakde, Jesse Silliman, and Jiuya Wang.
\newblock The {B}rumer-{S}tark conjecture over $\mathbb{Z}$.
\newblock {\em {P}reprint, arXiv:2310.16399}, 2023.

\bibitem{Sage/Vallieres:2022}
Sage DuBose and Daniel Valli\`eres.
\newblock On {$\mathbb{Z}^{d}_{\ell}$}-towers of graphs.
\newblock {\em Algebr. Comb.}, 6(5):1331--1346, 2023.

\bibitem{Gambheera/Vallieres:2024}
Rusiru Gambheera and Daniel Valli\`eres.
\newblock Iwasawa theory for branched {$\mathbb{Z}_p$}-towers of finite graphs.
\newblock {\em Doc. Math.}, 29(6):1435--1468, 2024.

\bibitem{Ghosh/Ray:2025}
Sohan Ghosh and Anwesh Ray.
\newblock On the {I}wasawa theory of {C}ayley graphs.
\newblock {\em Res. Math. Sci.}, 12(1):Paper No. 2, 2025.

\bibitem{Gonet:2022}
Sophia~R. Gonet.
\newblock Iwasawa theory of {J}acobians of graphs.
\newblock {\em Algebr. Comb.}, 5(5):827--848, 2022.

\bibitem{HMSV:2024}
Kyle Hammer, Thomas~W. Mattman, Jonathan~W. Sands, and Daniel Valli\`{e}res.
\newblock The special value {$u=1$} of {A}rtin-{I}hara {$L$}-functions.
\newblock {\em Proc. Amer. Math. Soc.}, 152(2):501--514, 2024.

\bibitem{Hashimoto:1990}
Ki-ichiro Hashimoto.
\newblock On zeta and {$L$}-functions of finite graphs.
\newblock {\em Internat. J. Math.}, 1(4):381--396, 1990.

\bibitem{Herbrand:1932}
Jacques Herbrand.
\newblock Sur les classes des corps circulaires.
\newblock {\em Journal de Math\'{e}matiques Pures et Appliqu\'{e}es},
  11:417--441, 1932.

\bibitem{Ihara:1966}
Yasutaka Ihara.
\newblock On discrete subgroups of the two by two projective linear group over
  {$\mathfrak{p}$}-adic fields.
\newblock {\em J. Math. Soc. Japan}, 18:219--235, 1966.

\bibitem{Kataoka:2024}
Takenori Kataoka.
\newblock Fitting ideals of {J}acobian groups of graphs.
\newblock {\em Algebr. Comb.}, 7(3):597--625, 2024.

\bibitem{Kleine/Muller:2023}
S\"{o}ren Kleine and Katharina M\"{u}ller.
\newblock On the growth of the {J}acobian in {$\mathbb{Z}_{p}^{l}$}-voltage
  covers of graphs.
\newblock {\em Algebr. Comb.}, 7(4):1011--1038, 2024.

\bibitem{Kleine/Muller:2025}
S\"{o}ren Kleine and Katharina M\"{u}ller.
\newblock On the non-commutative {I}wasawa main conjecture for voltage covers
  of graphs.
\newblock {\em {T}o appear in Israel Journal of Mathematics}, 2025.

\bibitem{Kundu/Muller:2024}
Debanjana Kundu and Katharina M\"{u}ller.
\newblock Iwasawa theory of graphs and their duals.
\newblock {\em {P}reprint, arXiv:2410.11704}, 2024.

\bibitem{Lei/Muller2}
Antonio Lei and Katharina M\"{u}ller.
\newblock On ordinary isogeny graphs with level structures.
\newblock {\em Expo. Math.}, 42(5):Paper No. 125589, 30, 2024.

\bibitem{Lei/Muller1}
Antonio Lei and Katharina M\"{u}ller.
\newblock On the zeta functions of supersingular isogeny graphs and modular
  curves.
\newblock {\em Arch. Math. (Basel)}, 122(3):285--294, 2024.

\bibitem{Lei/Muller3}
Antonio Lei and Katharina M\"{u}ller.
\newblock On towers of isogeny graphs with full level structures.
\newblock {\em Res. Math. Sci.}, 12(1):Paper No. 4, 29, 2025.

\bibitem{lei2022non}
Antonio Lei and Daniel Valli\`eres.
\newblock The non-{$\ell$}-part of the number of spanning trees in abelian
  {$\ell$}-towers of multigraphs.
\newblock {\em Res. Number Theory}, 9(1):Paper No. 18, 16, 2023.

\bibitem{Mazur:1977}
Barry Mazur.
\newblock Modular curves and the {E}isenstein ideal.
\newblock {\em Inst. Hautes \'{E}tudes Sci. Publ. Math.}, (47):33--186 (1978),
  1977.
\newblock With an appendix by Mazur and M. Rapoport.

\bibitem{Mazur-Wiles}
Barry Mazur and Andrew Wiles.
\newblock Class fields of abelian extensions of {${\bf Q}$}.
\newblock {\em Invent. Math.}, 76(2):179--330, 1984.

\bibitem{mcgownvallieresIII}
Kevin McGown and Daniel Valli\`eres.
\newblock On abelian {$\ell$}-towers of multigraphs {III}.
\newblock {\em Ann. Math. Qu\'{e}.}, 48(1):1--19, 2024.

\bibitem{Mizuno:2025}
Kosuke Mizuno.
\newblock Spanning trees and their relations in {G}alois covers.
\newblock {\em {P}reprint, arXiv:2503.19641}, 2025.

\bibitem{Northcott:1976}
D.~G. Northcott.
\newblock {\em Finite free resolutions}.
\newblock Cambridge Tracts in Mathematics, No. 71. Cambridge University Press,
  Cambridge-New York-Melbourne, 1976.

\bibitem{Pengo/Vallieres:2025}
Riccardo Pengo and Daniel Valli\`{e}res.
\newblock Spanning trees in {$\mathbb{Z}$}-covers of a finite graph and
  {M}ahler measures.
\newblock {\em J. Aust. Math. Soc.}, 118(1):108--144, 2025.

\bibitem{Ray/Vallieres:2023}
Anwesh Ray and Daniel Valli{\`e}res.
\newblock An analogue of {K}ida's formula in graph theory.
\newblock {\em To appear in {P}ure and {A}pplied {M}athematics {Q}uarterly},
  2025.

\bibitem{Ribet:1976}
Kenneth~A. Ribet.
\newblock A modular construction of unramified {$p$}-extensions of {$Q(\mu
  _{p})$}.
\newblock {\em Invent. Math.}, 34(3):151--162, 1976.

\bibitem{Serre:1977}
Jean-Pierre Serre.
\newblock {\em Arbres, amalgames, {${\rm SL}_{2}$}}.
\newblock Ast\'{e}risque, No. 46. Soci\'{e}t\'{e} Math\'{e}matique de France,
  Paris, 1977.
\newblock Avec un sommaire anglais, R\'{e}dig\'{e} avec la collaboration de
  Hyman Bass.

\bibitem{Stark/Terras:1996}
Harold~M. Stark and Audrey~A. Terras.
\newblock Zeta functions of finite graphs and coverings.
\newblock {\em Adv. Math.}, 121(1):124--165, 1996.

\bibitem{Stark/Terras:2000}
Harold~M. Stark and Audrey~A. Terras.
\newblock Zeta functions of finite graphs and coverings. {II}.
\newblock {\em Adv. Math.}, 154(1):132--195, 2000.

\bibitem{Stark/Terras:2007}
Harold~M. Stark and Audrey~A. Terras.
\newblock Zeta functions of finite graphs and coverings. {III}.
\newblock {\em Adv. Math.}, 208(1):467--489, 2007.

\bibitem{Sunada:2013}
Toshikazu Sunada.
\newblock {\em Topological crystallography}, volume~6 of {\em Surveys and
  Tutorials in the Applied Mathematical Sciences}.
\newblock Springer, Tokyo, 2013.
\newblock With a view towards discrete geometric analysis.

\bibitem{Terras:2011}
Audrey Terras.
\newblock {\em Zeta functions of graphs}, volume 128 of {\em Cambridge Studies
  in Advanced Mathematics}.
\newblock Cambridge University Press, Cambridge, 2011.
\newblock A stroll through the garden.

\bibitem{SAGE}
{The Sage Developers}.
\newblock {\em {S}ageMath, the {S}age {M}athematics {S}oftware {S}ystem
  ({V}ersion 9.6)}, 2022.

\bibitem{Washington:1997}
Lawrence~C. Washington.
\newblock {\em Introduction to cyclotomic fields}, volume~83 of {\em Graduate
  Texts in Mathematics}.
\newblock Springer-Verlag, New York, second edition, 1997.

\bibitem{Webb:2016}
Peter Webb.
\newblock {\em A course in finite group representation theory}, volume 161 of
  {\em Cambridge Studies in Advanced Mathematics}.
\newblock Cambridge University Press, Cambridge, 2016.

\bibitem{Wiles:1980}
Andrew Wiles.
\newblock Modular curves and the class group of {${\bf Q}(\zeta _{p})$}.
\newblock {\em Invent. Math.}, 58(1):1--35, 1980.

\end{thebibliography}
\end{document}